%% file: PHFol-Potrie-Revised2.tex
\documentclass[a4paper,12pt]{amsart}

\usepackage{amssymb,amscd,amsthm,mathrsfs}
\usepackage[ansinew]{inputenc}
\usepackage[centertags]{amsmath}
\usepackage{epsfig,graphicx}

  \def\CC{{\mathbb C}} \def\DD{{\mathbb D}}

\def\QQ{{\mathbb Q}} \def\RR{{\mathbb R}}  \def\TT{{\mathbb T}}
   
 \def\ZZ{{\mathbb Z}}

    \def\cS{\mathcal{S}}
    
\def\cC{\mathcal{C}}   \def\cO{\mathcal{O}} \def\cU{\mathcal{U}}
  \def\cJ{\mathcal{J}}  
\def\cE{\mathcal{E}}    \def\cW{\mathcal{W}}
\def\cF{\mathcal{F}}  \def\cL{\mathcal{L}}

\newtheorem*{teo*}{Theorem}
\newtheorem*{prop*}{Proposition}
\newtheorem*{cor*}{Corollary}
\newtheorem*{goal*}{Goal}

\newtheorem*{teoA}{Theorem A}
\newtheorem*{teoB}{Theorem B}

\newtheorem{teo}{Theorem}[section]

\newtheorem{quest}{Question}
\newtheorem{cor}[teo]{Corollary}
\newtheorem*{af}{Claim}
\newtheorem{lema}[teo]{Lemma}
\newtheorem{prop}[teo]{Proposition}

\newcommand{\bi}{\begin{itemize}}
\newcommand{\ei}{\end{itemize}}

\theoremstyle{definition}
\newtheorem{defi}[teo]{Definition}
\theoremstyle{remark}
\newtheorem{obs}[teo]{Remark}

\newcommand{\demo}[1]{\vspace{.05in}{\sc\noindent Proof #1.}}
\newcommand{\dem}{\vspace{.05in}{\sc\noindent Proof.  }}
\newcommand{\lqqd}{\par\hfill {$\Box$} \vspace*{.05in}}
\newcommand{\finobs}{\par\hfill{$\diamondsuit$} \vspace*{.05in}}

\newcommand{\eps}{\varepsilon}

\newcommand{\en}{\subset}
\newcommand{\trans}{\pitchfork}

\DeclareMathOperator{\rango}{rank} 
 \DeclareMathOperator{\Vol}{Vol}
\DeclareMathOperator{\length}{length}
\DeclareMathOperator{\diam}{diam}

\author[R. Potrie]{Rafael Potrie}
\address{CMAT, Facultad de Ciencias, Universidad de la Rep\'ublica, Uruguay}
\urladdr{www.cmat.edu.uy/$\sim$rpotrie}
\email{rpotrie@cmat.edu.uy}

\title[Partial hyperbolicity and foliations in $\TT^3$]{Partial hyperbolicity and foliations in $\TT^3$}
\thanks{The author was partially supported by ANR Blanc DynNonHyp BLAN08-2$\_$313375, ANII's doctoral scholarship, FCE-3-2011-1-6749 and Balzan's research project of J.Palis.}

\begin{document}
\begin{abstract}
We prove that dynamical coherence is an open and closed property
in the space of partially hyperbolic diffeomorphisms of $\TT^3$
isotopic to Anosov. Moreover, we prove that strong partially
hyperbolic diffeomorphisms of $\TT^3$ are either dynamically
coherent or have an invariant two-dimensional torus which is
either contracting or repelling. We develop for this end some
general results on codimension one foliations which may be of
independent interest.

\bigskip



\end{abstract}

\maketitle


\section{Introduction}

\subsubsection{} It is well known that robust dynamical properties
have implications on the existence of $Df$-invariant geometric
structures (see for example \cite{ManheECL,DPU,BDP}). In
particular, partial hyperbolicity plays a fundamental role in the
study of robust transitivity and stable ergodicity (see \cite{BDV}
chapters 7 and 8).

In this paper, we intend to develop some topological implications
of having $Df$-invariant geometric structures in the case the
manifold is $\TT^3$ and we are interested in giving conditions
under which the $Df$-invariant bundles are integrable.

As in the hyperbolic (Anosov) case, invariant  foliations play a
substantial role in the understanding of the dynamics of partially
hyperbolic systems (even if many important results manage to avoid
the use of the existence of such foliations, see for example
\cite{BuW, BI, Chab}). A still (essentially) up to date account on
the results on invariant foliations can be found in the seminal treaties \cite{BP,HPS} (see also
\cite{Berger}).

Two kinds of partially hyperbolic systems of $\TT^3$ are treated
in this paper: We shall call \emph{partially hyperbolic} to
diffeomorphisms whose tangent bundle splits into two
$Df$-invariant bundles $E^{cs} \oplus E^u$ where $E^u$ is one
dimensional and uniformly expanded and $E^{cs}$ is two-dimensional
and its behavior dominated by the expansion in $E^u$. A
diffeomorphism will be said to be \emph{strongly partially
hyperbolic} if the tangent bundle splits into 3-one dimensional
$Df$-invariant bundles $E^s \oplus E^c \oplus E^u$ such that the
extreme bundles are respectively uniformly contracted and expanded
by $Df$ and the center one has an intermediate behavior. (See
Section \ref{SectionDefiniciones} for precise definitions as well
as related ones.)

We say that a (strongly) partially hyperbolic diffeomorphism $f$ is \emph{dynamically coherent} if it admits
$f$-invariant foliations tangent to its $Df$-invariant distributions with the possible exception of $E^s \oplus E^u$ in the strong partially hyperbolic  case.

\subsubsection{} In principle, dynamical coherence may not be neither
open nor closed property among partially hyperbolic dynamics.
There are some hypothesis that guarantee openness which are not
known to hold in general and are usually hard to verify (see
\cite{HPS,Berger}). We prove here the following result (see
Section \ref{SectionResultados} for complete statement of the
results in this paper):

\begin{teo*}
Dynamical coherence is an open and closed property among partially
hyperbolic diffeomorphisms of $\TT^3$ isotopic to linear Anosov
automorphisms. Moreover, the unique obstruction for a strong partially
hyperbolic diffeomorphism of $\TT^3$ (i.e. with splitting of the
form $T\TT^3=E^s\oplus E^c \oplus E^u$) to be dynamically coherent
is the existence of a contracting or repelling periodic
two-dimensional torus.
\end{teo*}

The second part of this result responds in the affirmative for
$\TT^3$ to a conjecture made by M.A. Rodriguez Hertz, F.Rodriguez
Hertz and R. Ures in general 3 dimensional manifolds (see
\cite{HHU}). It is important to remark that they have constructed
examples showing that dynamical coherence does not hold for every
strong partially hyperbolic diffeomorphism.  As a consequence of
our results, we obtain:

\begin{cor*}
Let $f: \TT^3 \to \TT^3$ be a strong partially hyperbolic
diffeomorphism such that $\Omega(f)=\TT^3$ then $f$ is dynamically
coherent. \end{cor*}

\subsubsection{} Recent results by Brin, Burago and Ivanov (see
\cite{BBI,BI}) allow one to obtain certain topological
descriptions from the existence of invariant foliations.

We remark that recently, \cite{BBI2} have shown that absolute
strong partially hyperbolic diffeomorphisms of $\TT^3$ are
dynamically coherent using a criterium of Brin (\cite{Brin}) which
relies in this stronger version of partial hyperborbolicity in an
essential way. This has been used by Hammerlindl in
\cite{Hammerlindl} to obtain leaf-conjugacy results for this kind
of systems. The proof of our results owes a lot to both of these
results.

Absolute partial hyperbolicity covers many important classes of
examples and this makes its study very important. However, the
definition of absolute domination is somewhat artificial and does
not capture the results about robust dynamical properties (such as
the results of \cite{DPU,BDP}). Absolute domination can be
compared with a pinching condition on the spectrum of the
differential much like the conditions used by Brin and Manning to
classify Anosov systems (\cite{BrinManning}). It is worth to
remark that our results are the first (non perturbative) dynamical
coherence results under the \emph{pointwise} definition.

%
%
%
%

In a forthcoming paper (\cite{HP}), with A. Hammerlindl we use the
techniques here as well as the ones developed in
\cite{Hammerlindl} in order to give a classification result for
strongly partially hyperbolic diffeomorphisms of $\TT^3$.

\smallskip
{\small {\bf Acknowledgements: }\textit{This paper is part of my
thesis (\cite{PotTesis}) under S. Crovisier and M. Sambarino, I
would like to express my gratitude to their invaluable patience,
time and support. Sylvain devoted a great deal of time in
correcting the first draft. This work benefited from conversations
with P.Berger, C. Bonatti, J.Brum, A. Hammerlindl, P. Lessa and
A.Wilkinson. The paper owes a lot to the referee who made the
unrewarding task of reading carefully the paper (many times) and
making a lot of important suggestions to improve the presentation
and the paper itself, I am very grateful for that.}}

\section{Definitions and some known results}\label{SectionDefiniciones}

\subsection{Partial hyperbolicity}\label{SectionBundlesPHandVH}
\subsubsection{} Given a $C^1$-diffeomorphism $f: M \to M$ and a $Df$-invariant
subbundle $E\en TM$ we say that $E$ is \emph{uniformly
contracting} (resp. \emph{uniformly expanding}) if there exists
$n_0>0$ (resp. $n_0<0$) such that:

$$ \|Df^{n_0}|_{E}\| < \frac 1 2 $$

Given two $Df$-invariant subbundles $E,F$ of $TM$, we shall say
that $F$ \emph{dominates} $E$ if there exists $n_0>0$ such that
for every $x\in M$ and any pair of unit vectors $v_E \in E(x)$ and
$v_F \in F(x)$ we have that:

$$ \|Df^{n_0} v_E \| < \frac 1 2 \|Df^{n_0} v_F\| $$

It is important to remark here that this notion of domination is
weaker that the one appearing in other literature. Sometimes this
concept is called \emph{pointwise} (or \emph{relative})
\emph{domination} in contraposition to \emph{absolute domination}
(see \cite{BP,HPS}).

We say that $F$ \emph{absolutely dominates} $E$ if there exists
$n_0>0$ and $\lambda>0$ such that for every $x\in M$ and any pair
of unit vectors $v_E \in E(x)$ and $v_F \in F(x)$ we have that:

 $$ \|Df^{n_0} v_E \| < \lambda < \|Df^{n_0} v_F\| $$

\subsubsection{}\label{SectionDefsVHPHSPH}
Consider a diffeomorphism $f: M \to M$ such that $TM = E_1 \oplus
\ldots \oplus E_k$ is a $Df$-invariant splitting of the tangent
bundle on $k\geq 2$ non trivial invariant subbundles such that
$E_i$ dominates $E_j$ if $i>j$. Following \cite{BDP} (see also
\cite{BDV} appendix B) we say that:

\begin{itemize}
\item[-] $f$ is \emph{partially hyperbolic} if either $E_1$ is
uniformly contracting or $E_k$ is uniformly expanding. \item[-]
$f$ is \emph{strongly partially hyperbolic} if both $E_1$ is
uniformly contracting and $E_k$ is uniformly expanding.
\end{itemize}

We add the word \emph{absolutely} before this concepts when the
domination involved is of absolute nature.

These properties are $C^1$-robust (see \cite{BDV} appendix B).

%

In dimension $3$ we may have the following forms of partial
hyperbolicity: $TM= E^s \oplus E^c \oplus E^u$ (which will be the
strongly partially hyperbolic case), $TM =E^s \oplus E^{cu}$ and $TM
= E^{cs} \oplus E^u$ (the partially hyperbolic case). For
simplicity and by the symmetry of the problem, we shall focus only
on the cases $TM=E^{cs}\oplus E^u$ and $TM =E^s \oplus E^c \oplus
E^u$. 

\subsection{Invariant and almost invariant foliations}\label{SectionFoliacionesInvariantesYCasi}
\subsubsection{}\label{SectionFolInestable} Along this paper, we shall understand by \emph{foliation} a
$C^0$-foliation with $C^1$-leaves which is tangent to a continuous
distribution. For a foliation $\cF$ in a manifold $M$, we shall
denote as $\cF(x)$ to the leaf of $\cF$ containing $x$.

It is a well known result that there always exist an invariant
foliation tangent to $E^u$ and $E^s$ (\cite{BP,HPS}). We will denote as $\cF^u$ to the \emph{unstable foliation} and $\cF^s$ the \emph{stable foliation}.

%

\subsubsection{}\label{SectionCoherenciaADC} We say that a partially hyperbolic diffeomorphism $f$ with
splitting $TM=E^{cs} \oplus E^u$ is \emph{dynamically coherent} if
there exists an $f$-invariant foliation $\cF^{cs}$ tangent to
$E^{cs}$.

When $f$ is strongly partially hyperbolic, we say that $f$ is
\emph{dynamically coherent} when both $E^s \oplus E^c$ and $E^c
\oplus E^u$ are tangent to $f$-invariant foliations $\cF^{cs}$ and
$\cF^{cu}$ respectively. See \cite{BuW2} for more discussion: in
particular, this implies the existence of a $f$-invariant
foliation $\cF^c$ tangent to $E^c$ (obtained as the intersection
of the foliations).

All known examples in dimension 3 verify the following property
which is clearly $C^1$-open (we show in Proposition
\ref{PropositionADCCerrado} that the property is also closed):

\begin{defi}[Almost dynamical coherence]\label{DefADC}
We say that $f: M \to M$ partially hyperbolic of the form $TM=
E^{cs} \oplus E^u$ is \emph{almost dynamically coherent} if there
exists a foliation $\cF$ transverse to the direction $E^u$.
\end{defi}

The introduction of this property was motivated by \cite{BI} where
it is shown to be verified by strong partially hyperbolic systems
in dimension 3.

\subsection{Isotopy class of a diffeomorphism of $\TT^3$}\label{SubSectionClaseIsotopia}

Let $f: \TT^3 \to \TT^3$ be a $C^1$-diffeomorphism and $\tilde f:
\RR^3 \to \RR^3$ a lift to the universal cover.

One can define $f_\ast$ as the matrix given by the automorphism
$(\tilde f(\cdot) -\tilde f(0)) : \ZZ^3 \to \ZZ^3$. This matrix
does not depend on the chosen lift of $f$. It is direct to check
that $f_\ast$ is represented by a matrix of integer coeficients
and determinant equal to $1$ or $-1$ (since one can see that
$(f_\ast)^{-1}= (f^{-1})_\ast$).

We say that $f$ is \emph{isotopic to Anosov}\footnote{The name has
to do with the fact that in this context $f$ is always isotopic to
$f_\ast$. However we shall not use this fact but we keep this
nomenclature to be coherent with the literature.} if the matrix
$f_\ast$ has no eigenvalues of modulus $1$.

It is well known that $f_\ast$ also represents the action of $f$
on the fundamental group of $\TT^3$ (which coincides with the
first homology group since it is abelian). We will sometimes view $f_\ast$ as a diffeomorphism of $\RR^3$ or $\TT^3$.

By compactness, we have that there exists $K_0>0$ such that for every
$x\in \RR^3$:

$$ d(\tilde f(x), f_\ast(x)) < K_0 $$


\section{Precise statement of results and organization of the paper}\label{SectionResultados}

\subsection{Statement of results}

\subsubsection{} We are now in conditions to state precisely our
results. In both of them we shall only assume that the diffeomorphisms are of class $C^1$, no further regularity is required.

Given two transversal foliations $\cF_1$ and $\cF_2$ in a simply
connected manifold $\tilde M$ we say that they have \emph{global
product structure} if for every $x,y \in \tilde M$ we have that
the intersection between $\cF_1(x)$ and $\cF_2(y)$ has exactly one
point.

\begin{teoA}
Let $f: \TT^3 \to \TT^3$ an almost dynamically coherent partially
hyperbolic diffeomorphism with splitting of the form $T\TT^3 =
E^{cs} \oplus E^u$ and with $\dim E^u=1$. Assume that $f$ is
isotopic to Anosov, then:
\begin{itemize}
\item[-] $f$ is (robustly) dynamically coherent and has a unique
$f$-invariant foliation $\cF^{cs}$ tangent to $E^{cs}$. \item[-]
The lifts of the foliations $\cF^{cs}$ and $\cF^u$ to the universal cover have global product structure. \item[-] If $f_\ast$ has two eigenvalues of modulus larger than $1$ then they must be real and different.
\end{itemize}
\end{teoA}

In section \ref{Section-GPSCuantitativo} we obtain some general
results on global product structure for codimension one foliations
which may be of independent interest.

As a consequence of the fact that almost dynamical coherence is an
open and closed property we obtain:

\begin{cor}
Dynamical coherence is an open and closed property among partially
hyperbolic diffeomorphisms of $\TT^3$ isotopic to Anosov.
\end{cor}


\subsubsection{} In the strong partially hyperbolic case we are able
to give a stronger result partly based on results of \cite{BI}
showing that strongly partially hyperbolic diffeomorphisms of
$3$-manifolds are always almost dynamically coherent.

\begin{teoB}
Let $f: \TT^3 \to \TT^3$ be a strong partially hyperbolic
diffeomorphism, then:
\begin{itemize}
\item[-] Either there exists a unique $f$-invariant foliation
$\cF^{cs}$ tangent to $E^s \oplus E^c$ or, \item[-] There exists a
periodic two-dimensional torus $T$ tangent to $E^s \oplus E^c$
which is repelling.
\end{itemize}
\end{teoB}

\begin{obs}\label{RemarkTorosAnosov}
Indeed, it is not hard to show that in the case there is a
repelling torus, it must be an \emph{Anosov torus} as defined in
\cite{HHU3} (this follows from Proposition 2.4 of \cite{BBI}). In
the example of \cite{HHU} it is shown that the second possibility
is not empty. \finobs
\end{obs}

A diffeomorphism $f$ is \emph{chain-recurrent} if there is no open
set $U$ such that $f(\overline{U})\en U$ (see \cite{Chab} Chapter
1). In particular, if $\Omega(f)=\TT^3$ then
$f$ is chain-recurrent.

\begin{cor} Let $f: \TT^3 \to \TT^3$ be a chain-recurrent strongly partially hyperbolic diffeomorphism.
Then, $f$ is dynamically coherent.
\end{cor}

%
%


\subsection{Organization of the paper}
This paper is organized as follows: In section
\ref{SectionPreliminares} we introduce some new concepts, some
known results and prove preliminary results which will be used
later. In section \ref{SectionReeblesT3} we characterize Reebless
foliations on $\TT^3$ by following the ideas introduced in
\cite{BBI2}. Later, in section \ref{Section-GPSCuantitativo} we
give conditions under which the lift of some codimension one
foliations to the universal cover have global product structure.
In particular, we obtain a quantitative version of a known result
which we believe may be of independent interest. Finally, in section \ref{Section-PHAnosov} we prove
Theorem A and in section \ref{Section-TeoremaB} we prove Theorem
B.

In the appendix we give a proof of dynamical coherence for strong
partially hyperbolic diffeomorphisms isotopic to Anosov. This
proof is independent from Theorem A and easier. Other consequences
(relevant to \cite{HP}) are also obtained using some results from
the paper.

If the reader is interested in Theorem A only, she may skip
subsections \ref{SubSection-Branching},
\ref{SubSection-ClaseIsotopPH} and subsection
\ref{SubSection-FurtherPropertiesReeblesT3}.

If instead the reader is interested in Theorem B only, he may skip
subsection \ref{SubSection-ADCopenClosed}, section
\ref{Section-GPSCuantitativo} and read Appendix \ref{Apendice1}
instead of section \ref{Section-PHAnosov} (after having read
section \ref{Section-TeoremaB}).

The author suggest that the reader with some intuition on
foliations should skip sections 5 and 6 and come back to them when
having understood the core of the proof.


\section{Preliminary results and definitions}\label{SectionPreliminares}

In this paper, we will be mainly concerned with foliations and
diffeomorphisms of $\TT^3$. We shall fix the usual euclidean
metric as the distance in $\RR^3$, and this will induce a metric
on $\TT^3= \RR^3/_{\ZZ^3}$ via $p: \RR^3 \to \TT^3$ the canonical
quotient by translations of $\ZZ^3$.

We shall always denote as $B_\eps(\cC)$ to the $\eps$-neighborhood
of a set $\cC$ with the metric we have just defined. Usually
$\tilde X$ will denote the lift of $X$ to the universal cover
(whatever $X$ is).

\subsection{Some generalities on foliations}\label{SubSection-FoliacionesGral}

\subsubsection{}

For codimension one foliations, the fact that leaves cannot
intersect in more than one point is usually proved by using the
following result due to Haefliger (see \cite{Solodov,Hector} for
the $C^0$-version):

\begin{prop}[Haefliger's argument]\label{Proposicion-ArgumentoHaefliger}
Let $\cF$ be a codimension one foliation of a simply connected
(not necessarily compact) manifold $\tilde M$ such that every leaf
of $\cF$ is simply connected and let $\cF^\perp$ be a transverse
foliation, then, no leaf of $\cF^\perp$ can intersect a leaf of
$\cF$ twice.
\end{prop}

\subsubsection{}
We say that a foliation $\cF$ is \emph{without holonomy} if each
leaf of $\cF$ has trivial holonomy. We will not define this, but
we give two conditions under which a codimension one foliation
$\cF$ of a compact manifold is without holonomy:

\begin{itemize}
\item[-] If every leaf of $\cF$ is simply connected. \item[-] If
every leaf of $\cF$ is closed and every two leaves are homotopic.
\end{itemize}

Reeb's stability theorem (see \cite{Hector,CandelConlon}) has the
following consequences which we will use later:

\begin{prop}\label{Proposition-ReebStability}
Let $\cF$ be a codimension one foliation of a $3$-dimensional
manifold $M$:
\begin{itemize}
\item[(i)] If there is a leaf homeomorphic to $S^2$ then $M$ is
finitely covered by $S^2 \times S^1$. \item[(ii)] If a foliation
$\cF$ of $\TT^3$ has a dense set of leaves homeomorphic to
two-dimensional torus which are all homotopic, then, every leaf of
$\cF$ is a torus homotopic to the ones in the dense set. In
particular, $\cF$ is without holonomy.
\end{itemize}
\end{prop}

\subsubsection{} In \cite{BW} the following criterium for obtaining a
foliation out of a dense subset of leaves was given:

\begin{prop}[\cite{BW} Proposition 1.6 and Remark 1.10]\label{Proposition-BWProp16}
Let $E$ be a continuous codimension one distribution on a manifold
$M$ and $S$ a (possibly non connected) injectively immersed
submanifold everywhere tangent to $E$ which contains a family of
disks of fixed radius and whose set of midpoints is dense in $M$.
Then, there exists a foliation $\cF$ tangent to $E$ which contains
$S$ in its leaves.
\end{prop}

\subsection{Almost dynamical coherence is an open and closed property}\label{SubSection-ADCopenClosed}

\subsubsection{} Almost dynamical coherence is not a very strong
requirement. From its definition and the continuous variation of
the unstable bundle under perturbations it is clear that it is an
open property. With the basic facts on domination we can show:

\begin{prop}\label{PropositionADCCerrado}
Let $\{f_n\}$ a sequence of almost dynamically coherent partially
hyperbolic diffeomorphisms converging in the $C^1$-topology to a
partially hyperbolic diffeomorphism $f$. Then $f$ is almost
dynamically coherent.
\end{prop}

\dem Let us call $E^{cs}_n \oplus E^u_n$ to the splitting of $f_n$
and $E^{cs}\oplus E^u$ to the splitting of $f$. 
%

Consider $f_n$ such that the angle between $E^{cs}_n$ and
$E^{u}$ is everywhere larger than $\alpha/2$. Let $\cF_n$ be the
foliation transverse to $E^u_n$ given by the fact that $f_n$ is
almost dynamically coherent.


We can choose $m>0$ sufficiently large such that $f_n^{-m}(\cF_n)$
is tangent to a small cone (of angle less than $\alpha/2$) around
$E^{cs}_n$.

This implies that $f_n^{-m}(\cF_n)$ is a foliation transverse to
$E^u$ and this gives that $f$ is almost dynamically coherent as
desired.
\lqqd
%

\begin{obs} It is a mayor problem to determine whether partially hyperbolic
diffeomorphisms in $3$-dimensional manifolds are almost
dynamically coherent.
\end{obs}

%

\subsection{Consequences of Novikov's Theorem}\label{SectionNovikovBBI}

\subsubsection{} We shall state some results motivated by the work of Brin, Burago
and Ivanov (\cite{BBI,BI,BBI2}) on strong partial hyperbolicity.
They made the beautiful remark that a foliation transverse to the
unstable foliation cannot have Reeb components since that would
imply the existence of a closed unstable curve which is
impossible. We shall not define Reeb component here, but we refer
the reader to \cite{CandelConlon} chapter 9 for information on
this concepts. Reeb components in invariant foliations for
partially hyperbolic diffeomorphisms where also considered in
\cite{DPU} (Theorem H) but in the context of robust transitivity
and with different arguments.


\subsubsection{} The remark of Brin, Burago and Ivanov can be coupled with the
$C^0$ version of the classical Novikov's theorem (see
\cite{CandelConlon} Theorems 9.1.3 and 9.1.4).  We collect here
some of the consequences as obtained in \cite{BI,BBI2} (many of
the statements here also hold for general 3-dimensional manifolds,
see for example \cite{Pw} for related results):

\begin{teo}\label{CorolarioConsecuenciasReeb}
Let $f$ be a partially hyperbolic diffeomorphism of $\TT^3$ of the
form $T\TT^3=E^{cs}\oplus E^u$ ($\dim E^{cs}=2$) which is almost
dynamically coherent with foliation $\cF$. Let $\tilde \cF$ and
$\tilde \cF^u$ the lifts of the foliations $\cF$ and the unstable
foliation $\cF^u$ to $\RR^3$. Then:
\begin{itemize}
\item[(i)] For every $x\in \RR^3$ we have that $\tilde \cF(x) \cap
\tilde \cF^u(x) = \{x\}$. \item[(ii)] The leafs of $\tilde \cF$
are properly embedded complete surfaces in $\RR^3$. In fact there
exists $\delta>0$ such that every euclidean ball $U$ of radius
$\delta$ can be covered by a continuous coordinate chart such that
the intersection of every leaf $S$ of $\tilde \cF$ with $U$ is
either empty or represented as the graph of a function $h_S: \RR^2
\to \RR$ in those coordinates. \item[(iii)] Each closed leaf of
$\cF$ is an incompressible two dimensional torus (i.e. such that
the inclusion induces an injection of fundamental groups).
\item[(iv)] For every $\delta>0$, there exists a constant
$C_\delta$ such that if $J$ is a segment of $\tilde \cF^u$ then
$\Vol(B_\delta(J)) > C_\delta \length(J)$.
\end{itemize}
\end{teo}

\dem If a foliation $\cF$ on a compact closed $3$-dimensional
manifold $M$ has a Reeb component, then, every one dimensional
foliation transverse to $\cF$ has a closed leaf (see \cite{BI}
Lemma 2.2). Since $\cF^u$ is one dimensional, transverse to $\cF$
and has no closed leafs, we obtain that $\cF$ cannot have Reeb
components. To prove (i) and (ii) one can assume that $\cF$ is
transversally oriented since this holds for a finite lift and the
statement is in the universal cover.

The proof of (i) is the same as the one of Lemma 2.3 of \cite{BI}.

Once (i) is proved, (ii) follows from the same argument as in
Lemma 3.2 in \cite{BBI2}. 

Part (iii) follows from the fact that if $S$ is a closed surface
in $\TT^3$ which is not a torus, then it is either a sphere or its
fundamental group cannot inject in $\TT^3$ (see
\cite{ClasificacionSuperficies} and notice that neither a group
with exponential growth nor the fundamental group of the Klein
bottle can inject in $\ZZ^3$).

Since $\cF$ has no Reeb components (and the same happens for any
finite lift) we obtain that if $S$ is a closed leaf of $\cF$ then
it must be a sphere or an incompressible torus. But $S$ cannot be
a sphere since in that case, the Reeb's stability theorem (see
Proposition \ref{Proposition-ReebStability}) would imply that all
the leafs of $\cF$ are spheres and that the foliated manifold is
finitely covered by $S^2 \times S^1$ which is not the case.

The proof of (iv) is as Lemma 3.3 of \cite{BBI2}. 
%

\lqqd

\subsubsection{} Recall from subsection \ref{SubSectionClaseIsotopia} that if
$f:\TT^3 \to \TT^3$ then $f_\ast \in GL(3,\ZZ)$ denotes its action
on the first homology group (which has determinant of modulus
$1$).

\begin{teo}[Brin-Burago-Ivanov \cite{BBI,BI}]\label{TeoremaBuragoIvanovPH}
Let $f: \TT^3 \to \TT^3$ an almost dynamically coherent partially
hyperbolic diffeomorphism. Then, $f_\ast$ has at least an
eigenvalue of absolute value larger than $1$ and at least one of
absolute value smaller than $1$.
\end{teo}

\subsection{Branching foliations and Burago-Ivanov's results}\label{SubSection-Branching}

We follow \cite{BI} section 4.

We define a \emph{surface} in a $3$-manifold $M$ to be a
$C^1$-immersion $\imath: U \to M$ of a connected smooth
2-dimensional manifold without boundary. The surface is said to be
\emph{complete} if it is complete with the metric induced in $U$
by the Riemannian metric of $M$ pulled back by the immersion
$\imath$.

Given a point $x$ in (the image of) a surface $\imath: U \to M$ we
have that there is a neighborhood $B$ of $x$ such that the
connected component $C$ containing $\imath^{-1}(x)$ of
$\imath^{-1}(B)$ verifies that $\imath(C)$ separates $B$. We say
that two surfaces $\imath_1:U_1 \to M, \imath_2:U_2\to M$
\emph{topologically cross} if there exists a point $x$ in (the
image of) $\imath_1$ and $\imath_2$ and a curve $\gamma$ in $U_2$
such that $\imath_2(\gamma)$ passes through $x$ and intersects
both connected components of a neighborhood of $x$ with the part
of the surface defined above removed. It is not hard to prove that
the definition is indeed symmetric (see \cite{BI} Section 4).

A \emph{branching foliation} on $M$ is a collection of complete
surfaces tangent to a given continuous 2-dimensional distribution
on $M$ such that:

\begin{itemize}
\item[-] Every point belongs to (the image of) at least one
surface. \item[-] There are no topological crossings between
surfaces of the collection. \item[-] The branching foliation is
complete in the following sense: If $x_k \to x$ and $L_k$ are
(images of) surfaces of the collection through the points $x_k$,
we have that $L_k$ converges in the $C^1$-topology to (the image
of) a surface $L$ of the collection through $x$.
\end{itemize}

We call the (image of the) surfaces \emph{leaves} of the branching
foliation. We will abuse notation and denote a branching foliation
as $\cF_{bran}$ and by $\cF_{bran}(x)$ the set of leaves which
contain $x$.

\begin{teo}[\cite{BI},Theorem 4.1 and Theorem 7.2]\label{TeoBuragoIvanov}
Let $f:M^3 \to M^3$ a strong partially hyperbolic diffeomorphism
with splitting $TM=E^s\oplus E^c \oplus E^u$ into non trivial
one-dimensional bundles. There exists branching
foliations\footnote{The fact that one can choose them complete in
the sense defined above is proved in Lemma 7.1 of \cite{BI}.}
$\cF^{cs}_{bran}$ and $\cF^{cu}_{bran}$ tangent to
$E^{cs}=E^s\oplus E^c$ and $E^{cu}=E^c \oplus E^u$ respectively
which are $f$-invariant. For every $\eps>0$ there exist foliations
$\cS_\eps$ and $\cU_\eps$ tangent to an $\eps$-cone around
$E^{cs}$ and $E^{cu}$ respectively. Moreover, there exist
continuous and surjective maps $h^{cs}_\eps$ and $h^{cu}_\eps$ at
$C^0$-distance smaller than $\eps$ from the identity which send
the leaves of $\cS_\eps$ and $\cU_\eps$ to leaves of
$\cF^{cs}_{bran}$ and $\cF^{cu}_{bran}$ respectively.
\end{teo}

Since $E^s$ and $E^u$ are uniquely integrable we get that leaves of $\cF^s$ and
$\cF^u$ are contained in the leaves of $\cF^{cs}_{bran}$ and
$\cF^{cu}_{bran}$ respectively.

\begin{obs}\label{Remark-CercaEnCubrimientoYcercaImagen}
Notice that the existence of the maps $h^{cs}_\eps$ and
$h^{cu}_\eps$ which are $\eps$-close to the identity implies that
when lifted to the universal cover, the leaves of $\cS_\eps$
(resp. $\cU_\eps$) remain at distance smaller than $\eps$ from
lifted leaves of $\cF^{cs}_{bran}$ (resp. $\cF^{cu}_{bran}$).
\finobs
\end{obs}

%
%
%

The proof of Theorem B relies on the following property of
branching foliations whose proof follows from Proposition
\ref{Proposition-BWProp16}.

\begin{prop}\label{PropositionBranchingSinBranchingEsFoliacion}
If every point of $M$ belongs to a unique leaf of the branching
foliation, then the branching foliation is a true foliation.
\end{prop}

%
%
%
%

\subsection{Diffeomorphisms isotopic to linear Anosov automorphisms}\label{SectionResultadosdeSemiconjugacionWalters}

\subsubsection{} Let $f: \TT^3 \to \TT^3$ be a diffeomorphism which is
isotopic to a linear Anosov automorphism $A: \TT^3 \to \TT^3$.

We shall denote as $A$ to both the diffeomorphism of $\TT^3$ and
to the hyperbolic matrix $A \in GL(3,\ZZ)$ (with determinant of
modulus $1$) which acts in $\RR^3$ and is the lift of the torus
diffeomorphism $A$ to the universal cover.


\subsubsection{}
%

We have the following well result (see \cite{W} or \cite{PotTesis}
Section 2.3):

\begin{prop}\label{PropExisteSemiconjugacion}
For $f$ as above, there exists $H: \RR^3 \to \RR^3$ continuous and
surjective such that $H \circ \tilde f = A \circ H$. Also, it is
verified that $H(x + \gamma) = H(x) + \gamma$ for every $x \in
\RR^3$ and $\gamma \in \ZZ^3$ so, there exists also $h: \TT^3 \to
\TT^3$ homotopic to the identity such that $h \circ f = A \circ
h$. In particular, we have that $d(H(x), x)< K_1$ for every $x\in
\RR^3$.
\end{prop}

%
%
%
%

It is well known and easy to show that $H(W^\sigma(x,\tilde f))
\en W^\sigma(H(x),A)$ with $\sigma= s,u$.

\subsubsection{}\label{SectIrreducAnosov}
We finish this section by reviewing some properties of linear
Anosov automorphisms on $\TT^3$ (see for example \cite{PotTesis}
Section 1.5 for proofs).

We say that a matrix $A \in GL(3,\ZZ)$ (with determinant of
modulus $1$) is \emph{irreducible} if and only if its
characteristic polynomial is irreducible in the field $\QQ$. It is
not hard to prove:

\begin{prop}\label{PropAnosovSonIrreducibles}
Every hyperbolic matrix $A \in GL(3,\ZZ)$ with determinant of
modulus $1$ is irreducible. In particular: it cannot have an
invariant linear two-dimensional torus; all its eigenvalues are
real, irrational and different and if there are two complex
conjugate eigenvalues they have irrational angle.
\end{prop}

\subsubsection{} Since the eigenspaces of the matrix $A$ are invariant
subspaces for $A$, we get that each eigenline or eigenplane has
totally irrational direction (by this we mean that its projection
to $\TT^3$ is simply connected and dense). Also, these are the
only invariant subspaces of $A$.

Since $A$ is hyperbolic and the product of eigenvalues is one, we
get that $A$ must have one or two eigenvalues with modulus smaller
than $1$. We say that $A$ has \emph{stable dimension} $1$ or $2$
depending on how many eigenvalues of modulus smaller than one it
has.

We call \emph{stable eigenvalues} (resp. \emph{unstable
eigenvalues}) to the eigenvalues of modulus smaller than one
(resp. larger than one). The subspace $E^s_A$ (resp $E^u_A$)
corresponds to the eigenspace associated to the stable (resp.
unstable) eigenvalues.


\begin{obs}\label{RemarkSubespaciosInvariantesAnosov}
For every $A \in GL(3,\ZZ)$ hyperbolic with determinant of modulus
$1$, we know exactly which are the invariant planes of $A$. If $A$
has complex eigenvalues, then, the only invariant plane is the
eigenspace associated to that pair of complex conjugate
eigenvalues. If $A$ has $3$ different real eigenvalues then there
are $3$ different invariant planes, one for each pair of
eigenvalues. All these planes are totally irrational. \finobs
\end{obs}

\subsection{Diffeomorphisms isotopic to linear partially hyperbolic maps of $\TT^3$}\label{SubSection-ClaseIsotopPH}

\subsubsection{} In the case $f_\ast$ is hyperbolic, we saw that each
eigendirection of $f_\ast$ projects into a immersed line which is
dense in $\TT^3$ (the same holds for each plane).

In the non hyperbolic partially hyperbolic case (i.e. when one
eigenvalue has modulus $1$ and the other two have modulus
different from $1$), we get that:

\begin{lema}\label{LemaPlanosInvariantes}
Let $A$ be a matrix in $GL(3,\ZZ)$ with eigenvalues $\lambda^s,
\lambda^c, \lambda^u$ verifying $0<|\lambda^s| < |\lambda^c|=1 <
|\lambda^u| = |\lambda^s|^{-1}$. Let $E^s_A, E^c_A, E^u_A$ be the
eigenspaces associated to $\lambda^s, \lambda^c$ and $\lambda^u$
respectively. We have that:
\begin{itemize}
\item[-] $E^c_A$ projects by $p$ into a closed circle where
$p:\RR^3 \to \TT^3$ is the covering projection. \item[-] The
eigenlines $E^{s}_A$ and $E^u_A$ project by $p$ into immersed
lines whose closure coincide with a two dimensional linear torus.
\end{itemize}
In particular, if $P$ is an $A$-invariant plane then its
projection to $\TT^3$ is either a torus or a dense cylinder.
\end{lema}

See \cite{PotTesis} Section 1.5 for a proof of this standard result.

\section{Reebless foliations of $\TT^3$}\label{SectionReeblesT3}

\subsection{Some preliminaries and statement of results}\label{SubSection-EnunciadosReeblesT3}

\subsubsection{} We consider a codimension one foliation $\cF$ of
$\TT^3$  and $\cF^\perp$ a one dimensional transversal foliation.
We shall assume throughout that $\cF$ has no Reeb
components\footnote{In case the reader is not familiar with Reeb
components, we remark that for the purposes of this paper, a
foliation will be \emph{Reebless} if the conclusions of Theorem
\ref{CorolarioConsecuenciasReeb} are satisfied.} and that
$\cF^\perp$ is oriented (we can always assume this by considering
a double cover).


Consider an orientation on $\tilde \cF^\perp$. Given $x\in \RR^3$ we get that $\tilde \cF^\perp(x) \setminus
\{x\}$ has two connected components which we call $\tilde
\cF^\perp_+(x)$ and $\tilde \cF^\perp_-(x)$ according to the
chosen orientation of $\tilde \cF^\perp$. We denote as
$F_+(x)$ and $F_-(x)$ to the connected components of $\RR^3 \setminus \tilde \cF(x)$ depending on whether they contain $\tilde
\cF^\perp_+(x)$ or $\tilde \cF^\perp_-(x)$.


Since covering transformations preserve the orientation and
$\tilde \cF$:

$$F_{\pm}(x) + \gamma = F_{\pm}(x+\gamma) \qquad \forall \gamma \in \ZZ^3 $$

\subsubsection{} For every $x\in \RR^3$, we consider the following subsets of $\ZZ^3$ seen as deck transformations:

$$\Gamma_+(x) = \{ \gamma \in \ZZ^3 \ : \ F_+(x) + \gamma \en F_+(x) \} $$

$$\Gamma_-(x) = \{ \gamma \in \ZZ^3 \ : \ F_-(x) + \gamma \en F_-(x) \} $$

We also consider $\Gamma(x) = \Gamma_+ (x) \cup \Gamma_-(x)$.


\begin{lema}\label{RemarkPosibilidadesdelosFmas}  The following properties hold:
\begin{itemize}
\item[(i)] If both $F_+(x)\cap F_+(y) \neq \emptyset$ and $F_-(x)
\cap F_-(y) \neq \emptyset$ then, either $F_+(x) \en F_+(y)$ and
$F_-(y) \en F_-(x)$ or $F_+(y) \en F_+(x)$ and $F_-(x) \en
F_-(y)$. In both of this cases we shall say that $F_+(x)$ and
$F_+(y)$ are \emph{nested} (similar with $F_-$). \item[(ii)] If
$F_+(x) \cap F_+(y)=\emptyset$ then $F_+(y) \en F_-(x)$ and
$F_+(x) \en F_-(y)$. A similar property holds if $F_-(x) \cap
F_-(y)=\emptyset$. \item[(iii)] In particular, $F_+(x) \en F_+(y)$
if and only if $F_-(y) \en F_-(x)$ (they are nested in both
cases).
\end{itemize}
\end{lema}

\dem We will only consider the case where $\tilde \cF(x) \neq
\tilde \cF(y)$ since otherwise the Lemma is trivially satisfied
(and case (ii) is not possible).

Assume that both $F_+(x) \cap F_+(y)$ and $F_-(x) \cap F_-(y)$ are
non-empty. Since $\tilde \cF(y)$ is connected and does not
intersect $\tilde \cF(x)$ we have that it is contained in either
$F_+(x)$ or $F_-(x)$. We can further assume that $\tilde \cF(y)
\en F_+(x)$ the other case being symmetric. In this case, we
deduce that $F_+(y) \en F_+(x)$: otherwise, we would have that
$F_-(x) \cap F_-(y) = \emptyset$. But this implies that $\tilde
\cF(x) \en F_-(y)$ and thus that $F_-(x) \en F_-(y)$ which
concludes the proof of (i).

To prove (ii) notice that if $F_+(x) \cap F_+(y) = \emptyset$ then
we have that $\tilde \cF(x) \en F_-(y)$ and $\tilde \cF(y) \en
F_-(x)$. This gives that both $F_+(x) \en F_-(y)$ and $F_+(y) \en
F_-(x)$ as desired.

Finally, if $F_+(x) \en F_+(y)$ we have that $F_-(x) \cap F_-(y)$
contains at least $F_-(y)$ so that (i) applies to give (iii).

\lqqd

\begin{figure}[ht]
\begin{center}
\input{fmasopciones.pstex_t}
\caption{\small{When $F_+(x)$ and $F_+(x)+\gamma$ are not
nested.}} \label{FiguraFmasOpciones}
\end{center}
\end{figure}
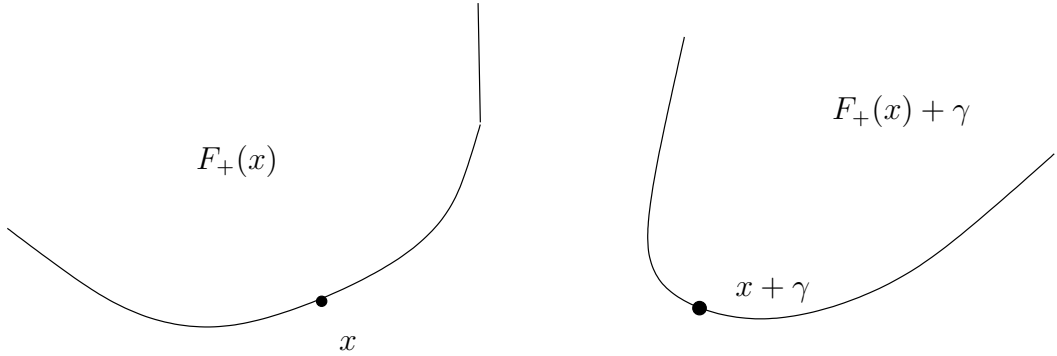


One can then prove (see Lemma 3.9 of \cite{BBI2}):

\begin{lema}\label{LemaGamaEsSubGrupoDeZ3} For every $x\in \RR^3$ we have that $\Gamma(x)$ is a subgroup of $\ZZ^3$.
\end{lema}

%
%
%
%
%

\subsubsection{} We close this subsection by stating the following theorem
which provides a kind of classification of Reebless foliations in
$\TT^3$ (see also \cite{Plante} where a similar result is obtained
for $C^2$-foliations of more general manifolds). The proof is
deferred to the next subsection.

We say that $\cF$ has a \emph{dead end component} if there exists
two (homotopic) torus leaves $T_1$ and $T_2$ of $\cF$ such that
there is no transversal that intersects both of them. When such a
component exists, we have that the leaves of any transversal
foliation must remain at bounded distance from some lift of $T_1$
and $T_2$.

\begin{teo}\label{PropObienHayTorosObienTodoEsLindo}
Let $\cF$ be a Reebless foliation of $\TT^3$. Then, there exists a
plane $P \en \RR^3$ and $R>0$ such that every leaf of $\tilde \cF$
lies in an $R$-neighborhood of a translate of $P$. Moreover, one
of the following condition holds:
\begin{itemize}
\item[(i)] Either for every $x \in \RR^3$ the $R$-neighborhood of
$\tilde \cF(x)$ contains $P+x$, or, \item[(ii)] $P$ projects into
a $2$-torus and there is a dead-end component of $\cF$.
\end{itemize}
\end{teo}

%
%

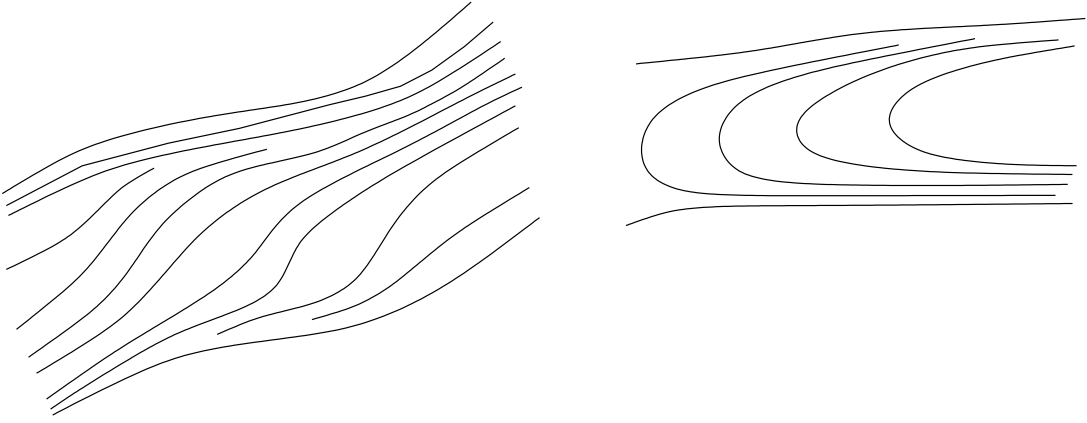
\begin{figure}[ht]
\begin{center}
\input{foliations.pstex_t}
\caption{\small{How the possibilities on $\tilde \cF$ look like.}}
\label{FiguraFmasOpciones}
\end{center}
\end{figure}

In case the foliation $\cF^\perp$ is not oriented, we get
essentially the same results.

\subsection{Proof of Theorem \ref{PropObienHayTorosObienTodoEsLindo}}

\subsubsection{}  We follow the proof in \cite{BBI2} with
some slight adaptations to get our statement.


For $x\in \RR^3$, define $G_+(x)=\bigcap_{\gamma \in \ZZ^3}\overline{F_+(x) +
\gamma}$ and $G_-(x)$ in a similar way.

First, assume that there exists $x_0$ such that $G_+(x_0) \neq
\emptyset$ (see Lemma 3.10 of \cite{BBI2}). The case where
$G_-(x_0)\neq \emptyset$ is symmetric. The idea is to prove that
in this case we will get option (ii) of the theorem.

%
%
%

One can check that the intersection $G_+(x_0)$ is a
$3$-dimensional submanifold of $\RR^3$ (modeled in the upper half
space) with boundary consisting of leaves of $\tilde \cF$ (since
the boundary components are always locally limits of local
leaves).
%

\begin{af}
If $G_+(x_0) \neq \emptyset$ then there exists plane $P$ and $R>0$
such that every leaf of $\tilde \cF(x)$ is contained in an
$R$-neighborhood of a translate of $P$ and whose projection to
$\TT^3$ is a two dimensional tori. Moreover, option (ii) of the
proposition holds.
\end{af}

\dem Since $G_+(x_0)$ is invariant under every integer
translation, we get that the boundary of $G_+(x_0)$ descends to a
closed surface in $\TT^3$ which is union of leaves of $\cF$.

By Novikov's Theorem (recall Theorem
\ref{CorolarioConsecuenciasReeb} (iii)) we get that those leaves
are two-dimensional torus whose fundamental group is injected by
the inclusion map.

This implies that they are at bounded distance of linear
embeddings of $\TT^2$ in $\TT^3$ and so their lifts lie within
bounded distance from a plane $P$ whose projection is a two
dimensional torus.

Since leafs of $\tilde \cF$ do not cross, the plane $P$ does not
depend on the boundary component. Moreover, every leaf of $\tilde
\cF(x)$ must be at bounded distance from a translate of $P$ since
every leaf of $\cF$ has a lift which lies within two given lifts
of some of the torus leafs.

Consider a point $x$ in the boundary of $G_+(x_0)$. We have that
$\tilde \cF(x)$ is at bounded distance from $P$ from the argument
above.

Moreover, each boundary component of $G_+(x_0)$ is positively
oriented in the direction which points inward to the interior of
$G_+(x_0)$ (recall that it is a compact $3$-manifold with
boundary).

We claim that $\eta_z$ remains between two translates of $P$ for
every $z\in \tilde \cF(x)$ and $\eta_z$ positive transversal to
$\tilde \cF$. Indeed, if this is not the case, then $\tilde
\eta_z$ would intersect other boundary component of $G_+(x_0)$
which is impossible since the boundary leafs of $G_+(x_0)$ point
inward to $G_+(x_0)$ (with the orientation of $\tilde \cF^\perp$).

Now, consider any point $z \in \RR^3$, and $\eta_z$ a positive
transversal which we assume does not remain at bounded distance
from $P$. Then it must intersect some translate of $\tilde
\cF(x)$, and the argument above applies. This is a contradiction.

The same argument works for negative transversals since once a
leaf enters $(G_+(x_0))^c$ it cannot reenter any of its
translates. We have proved that $p(G_+(x_0))$ contains a dead end
component. This proves the claim.
\finobs

Now, assume that (ii) does not hold, in particular $G_\pm(x)
=\emptyset$ for every $x$. Then, for every point $x$ we have that
$$\bigcup_{\gamma \in \ZZ^3} (F_+(x) + \gamma) =
\bigcup_{\gamma\in \ZZ^3}(F_-(x) + \gamma) = \RR^3$$

\noindent As in Lemma 3.11 of \cite{BBI2} we can prove:

\begin{af} We have that $\Gamma(x)=\ZZ^3$ for every $x\in \RR^3$.
\end{af}
%
%
%
%
%
%
%
%

Consider $\Gamma_0(x)= \Gamma_+(x) \cap \Gamma_-(x)$, the set of
translates which fix $\tilde \cF(x)$.

If $\rango(\Gamma_0(x))=3$, then $p^{-1}(p(\tilde \cF(x)))$
consists of finitely many translates of $\tilde \cF(x)$ which
implies that $p(\tilde \cF(x))$ is a closed surface of $\cF$. On
the other hand, the fundamental group of this closed surface
should be isomorphic to $\ZZ^3$ which is impossible since there
are no closed surfaces with such fundamental group
(\cite{ClasificacionSuperficies}). This implies that
$\rango(\Gamma_0(x))<3$ for every $x\in \RR^3$.

\begin{af}
For every $x\in \RR^3$ there exists a plane $P(x)$ and translates
$P_+(x)$ and $P_-(x)$ such that $F_+(x)$ lies in a half space
bounded by $P_+(x)$ and $F_-(x)$ lies in a half space bounded by
$P_-(x)$.
\end{af}

\dem The argument is the same as in Lemma 3.12 of \cite{BBI2} (and the
argument after that lemma).

\finobs

We have proved that for every $x \in \RR^3$ there exists a plane
$P(x)$ and translates $P_+(x)$ and $P_-(x)$ such that $F_\pm (x)$
lies in a half space bounded by $P_{\pm}(x)$. Let $R(x)$ be the
distance between $P_+(x)$ and $P_-(x)$, we have that $\tilde
\cF(x)$ lies at distance smaller than $R(x)$ from $P_+(x)$.

Now, we must prove that the $R(x)$-neighborhood of $\tilde \cF(x)$
contains $P_+(x)$. To do this, it is enough to show that the
projection from $\tilde \cF(x)$ to $P_+(x)$ by an orthogonal
vector to $P(x)$ is surjective. If this is not the case, then
there exists a segment joining $P_+(x)$ to $P_-(x)$ which does not
intersect $\tilde \cF(x)$. This contradicts the fact that every
curve from $F_-(x)$ to $F_+(x)$ must intersect $\tilde \cF(x)$.

Since the leaves of $\tilde \cF$ do not intersect, $P(x)$ cannot
depend on $x$. Since the foliation is invariant under integer
translations, we get (by compactness) that $R(x)$ can be chosen
uniformly bounded. This concludes the proof of Theorem
\ref{PropObienHayTorosObienTodoEsLindo}.

\lqqd

\begin{obs}\label{Remark-UnicidaddelPlanoP}
It is direct to show that for a given Reebless foliation $\cF$ of
$\TT^3$, the plane $P$ given by Theorem
\ref{PropObienHayTorosObienTodoEsLindo} is unique. 
\end{obs}


\subsection{Further properties of the foliations}\label{SubSection-FurtherPropertiesReeblesT3}

\subsubsection{}

Given a codimension one foliation $\tilde \cF$ of $\RR^3$ whose
leaves are homeomorphic to two dimensional planes, one defines the
\emph{leaf space} $\cL$ of $\tilde \cF$ by the quotient space of
$\RR^3$ with the equivalence relationship of being in the same
leaf of $\tilde \cF$. It is well known that in this
case\footnote{When all leaves are properly embedded copies of
$\RR^2$.}, we have that $\cL= \RR^3 /_{\tilde \cF}$ is a (possibly
non-Hausdorff) one-dimensional manifold (see \cite{CandelConlon}
Proposition II.D.1.1).

It is not hard to see that:

\begin{prop}\label{RemarkElCocienteDeLaFoliacionEsR}
Let $\cF$ be a Reebless foliation of $\TT^3$, if option (i) of
Theorem \ref{PropObienHayTorosObienTodoEsLindo} holds then the
leaf space $\cL= \RR^3 /_{\tilde \cF}$ is homeomorphic to $\RR$.
\end{prop}

\dem The space of leafs $\cL$ with the quotient topology has the
structure of a (possibly non-Hausdorff) one-dimensional manifold.
In fact, this follows directly from Proposition
\ref{Proposicion-ArgumentoHaefliger} which also implies that it is
simply connected as a one-dimensional manifold. To prove the
proposition is thus enough to show that it is Hausdorff.

We define an ordering in $\cL$ as follows

$$\tilde \cF(x) \geq \tilde \cF(y) \quad \text{if} \quad  F_+(x) \en F_+(y).$$

If option (i) of Theorem \ref{PropObienHayTorosObienTodoEsLindo}
holds, given $x,y$ we have that $F_+(x)\cap F_+(y) \neq \emptyset$
and $F_-(x) \cap F_-(y) \neq \emptyset$.

Then, Lemma \ref{RemarkPosibilidadesdelosFmas} (i) implies that
$F_+(x)$ and $F_+(y)$ are nested. In conclusion, we obtain that
the relationship we have defined is a total order.

Let $\tilde \cF(x)$ and $\tilde \cF(y)$ two different leaves of
$\tilde \cF$. We must show that they belong to disjoint open sets.

Without loss of generality, since it is a total order, we can
assume that $\tilde \cF(x) < \tilde \cF(y)$. This implies that
$F_+(y)$ is strictly contained in $F_+(y)$. On the other hand,
this implies that $F_-(y) \cap F_+(x) \neq \emptyset$, in
particular, there exists $z$ such that $\tilde \cF(x) < \tilde
\cF(z) < \tilde \cF(y)$.

Since the sets $F_+(z)$ and $F_-(z)$ are open and disjoint and we
have that $\tilde \cF(x) \en F_-(z)$ and $\tilde \cF(y) \in
F_+(z)$ we deduce that $\cL$ is Hausdorff as desired.

\lqqd

Since $\tilde \cF$ is invariant under deck transformations, we
obtain that we can consider the quotient action of $\ZZ^3 =
\pi_1(\TT^3)$ in $\cL$. For $[x]=\tilde \cF(x) \in \cL$ we get
that $\gamma \cdot [x] = [x+\gamma]$ for every $\gamma \in \ZZ^3$.

\subsubsection{}
Notice that all leaves of $\cF$ in $\TT^3$ are simply connected if
and only if $\pi_1(\TT^3)$ acts without fixed point in $\cL$. In a
similar fashion, existence of fixed points, or common fixed points
allows one to determine the topology of leaves of $\cF$ in
$\TT^3$.

In fact, we can prove:

\begin{prop}\label{Proposicion-SiPesToroHayHojaToro}
Let $\cF$ be a Reebless foliation of $\TT^3$. If the plane $P$
given by Theorem \ref{PropObienHayTorosObienTodoEsLindo} projects
into a two dimensional torus by $p$, then there is a leaf of $\cF$
homeomorphic to a two-dimensional torus.
\end{prop}

\dem Notice first that if option (ii) of Theorem
\ref{PropObienHayTorosObienTodoEsLindo} holds, the existence of a
torus leaf is contained in the statement of the theorem.

So, we can assume that option (i) holds. By considering a finite
index subgroup, we can further assume that the plane $P$ is
invariant under two of the generators of $\pi_1(\TT^3)\cong \ZZ^3$
which we denote as $\gamma_1$ and $\gamma_2$.

Since leaves of $\tilde \cF$ remain close in the Hausdorff
topology to the plane $P$ we deduce that the orbit of every point
$[x] \in \cL$ by the action of the elements $\gamma_1$ and
$\gamma_2$ is bounded.

Let $\gamma_3$ be the third generator: its orbit cannot be
bounded, otherwise translation by $\gamma_3$ would fix the plane
$P$. So, the quotient of $\cL$ by the action of $\gamma_3$ is a
circle. We can make the group generated by $\gamma_1$ and
$\gamma_2$ act on this circle and we obtain two commuting circle
homeomorphisms with zero rotation number. This implies they have a
common fixed point which in turn gives us the desired two-torus
leaf of $\cF$.

\lqqd

\subsubsection{}

Also, depending on the topology of the projection of the plane $P$
given by Theorem \ref{PropObienHayTorosObienTodoEsLindo} we can
obtain some properties on the topology of the leaves of $\cF$.

\begin{lema}\label{RemarkHojasCerradas}
Let $\cF$ be a Reebless foliation of $\TT^3$ and $P$ be the plane
given by Theorem \ref{PropObienHayTorosObienTodoEsLindo}.
\begin{itemize}
\item[(i)] Every closed curve in a leaf of $\cF$ is homotopic in
$\TT^3$ to a closed curve contained in $p(P)$. This implies in
particular that if $p(P)$ is simply connected, then all the leaves
of $\cF$ are also simply connected. \item[(ii)] If a leaf of $\cF$
is homeomorphic to a two dimensional torus, then, it is homotopic
to $p(P)$ (in particular, $p(P)$ is also a two dimensional torus).
\end{itemize}
\end{lema}

\dem To see (i), first notice that leafs are incompressible. Given
a closed curve $\gamma$ in a leaf of $\cF$ which is not
null-homotopic, we know that when lifted to the universal cover it
remains at bounded distance from a linear one-dimensional subspace
$L$. Since $\gamma$ is a circle, we get that $p(L)$ is a circle in
$\TT^3$. If the subspace $L$ is not contained in $P$ then it must
be transverse to it. This contradicts the fact that leaves of
$\cF$ remain at bounded distance from $P$.

To prove (ii), notice that a torus leaf $T$ which is
incompressible must remain close in the universal cover to a plane
$P_T$ which projects to a linear embedding of a $2$-dimensional
torus. From the proof of Theorem
\ref{PropObienHayTorosObienTodoEsLindo} and the fact that $\cF$ is
a foliation we get that $P_T=P$. See also the proof of Lemma 3.10
of \cite{BBI2}.

\lqqd

\subsubsection{}

We end this section by obtaining some results about branching
foliations we will use only in Section \ref{Section-TeoremaB}.

\begin{prop}\label{Proposition-SucesionDeTorosEnBranchingFol}
Let $\cF_{bran}$ be a branching foliation of $\TT^3$ and consider
a sequence of points $x_k$ such that there are leaves $L_k \in
\cF_{bran}(x_k)$ which are closed, incompressible and homotopic to
each other. If $x_k \to x$, then there is a leaf $L \in
\cF_{bran}(x)$ which is incompressible and homotopic to the leaves
$L_k$.
\end{prop}

\dem Recall that if $x_k \to x$ and we consider a sequence of
leaves through $x_k$ we get that the leaves converge to a leaf
through $x$.

Since the leaves of $\cF_{bran}$ are incompressible, the lifts of
every leaf in $\cF_{bran}(x_k)$ is homeomorphic to a plane.
Moreover, the fundamental group of each leaf must be $\ZZ^2$ and
the leaves must be homeomorphic to 2-torus, since it is the only
possibly incompressible surface in $\TT^3$.

Since all the leaves $\cF_{bran}(x_k)$ are homotopic, their lifts
are invariant under the same elements of $\pi_1(\TT^3)$. The limit
leaf must thus be also invariant under those elements. Notice that
it cannot be invariant under further elements of $\pi_1(\TT^3)$
since no surface has such fundamental group.

\lqqd

\section{Global product structure: Quantitative results}\label{Section-GPSCuantitativo}

\subsection{Statement of results}

Recall that given two transverse foliations $\cF_1$ and $\cF_2$ of
a manifold $M$ we say they admit a \emph{global product structure}
if given two points $x,y \in \tilde M$ the universal cover of $M$
we have that $\tilde \cF_1(x)$ and $\tilde \cF_2(y)$ intersect in
a unique point.
%


The $C^2$-version of the following theorem is due to Novikov
according to \cite{Hector}:

\begin{teo}[Theorem VIII.2.2.1 of \cite{Hector}]\label{Teorema-HectorHirsch}
Consider a codimension one foliation $\cF$ without holonomy of a
compact manifold $M$. Then, for every $\cF^\perp$ foliation
transverse to $\cF$ we have that $\cF$ and $\cF^\perp$ have global
product structure.
\end{teo}

Other than the case where there is a compact leaf without
holonomy, the other important case in which this result applies is
when every leaf of the foliation is simply connected.
Unfortunately, there will be some situations where we will be
needing to obtain global product structure but not having neither
all leaves of $\cF$ simply connected nor that the foliation lacks
of holonomy in all its leaves.

We will use instead the following quantitative version of the
previous result which does not imply it other than it the
situations we will be needing it. We hope this general result on
the existence of global product structure may find other
applications.

\begin{teo}\label{Teorema-EstructuraProductoGlobalMIO}
Let $M$ be a compact manifold and $\delta>0$. Consider a set of
generators of $\pi_1(M)$ and endow $\pi_1(M)$ with the word length
for generators. Then, there exists $K>0$ such that if $\cF$ is a
codimension one foliation and $\cF^\perp$ a transverse foliation
such that:
\begin{itemize}
\item[-] There is a local product structure of size $\delta$
between $\cF$ and $\cF^\perp$. \item[-] The leaves of $\tilde \cF$
are simply connected and no element of $\pi_1(M)$ of size less
than $K$ fixes a leaf of $\tilde \cF$. \item[-] The leaf space
$\cL = \tilde M /_{\tilde \cF}$ is homeomorphic to $\RR$. \item[-]
$\pi_1(M)$ is abelian.
\end{itemize}
Then, $\cF$ and $\cF^\perp$ admit a global product structure.
\end{teo}

We remark that it is possible to construct a foliation of $\TT^3$
by planes whose leaf space in the universal cover is homeomorphic
to $\RR$ and which has a transverse one dimensional foliation
without global product structure in the universal cover. This
implies that the hypothesis of having ``small holonomy''  in the
sense of the Theorem (given the size of the local product
structure, one obtains a size of deck transformations without
fixed leafs) is essential to deduce the result.


\subsection{Proof of Theorem \ref{Teorema-EstructuraProductoGlobalMIO}}

Proposition \ref{Proposicion-ArgumentoHaefliger} (Haefliger's argument) implies that  leaves of $\tilde \cF$ and $\tilde
\cF^\perp$ intersect in at most one point, so one must only show that the intersection is non-empty.
%
%

In $\cL = \tilde M /_{\tilde \cF}$ we can consider an ordering of
leafs (by using the ordering from $\RR$). We denote as $[x]$ to
the equivalence class in $\tilde M$ of the point $x$, which
coincides with $\tilde \cF(x)$.

The following condition will be the main ingredient for obtaining
a global product structure:

\begin{itemize}
\item[($\ast$)] For every $z_0\in \tilde M$ there exists $y^-$ and
$y^+ \in \tilde M$ verifying that $[y^-]< [z_0] < [y^+]$ and such
that for every $z_1, z_2\in \tilde M$ satisfying $[y^-]\leq [z_i]
\leq [y^+]$ ($i=1,2$) we have that $\tilde \cF^\perp(z_1) \cap
\tilde \cF(z_2) \neq \emptyset.$
\end{itemize}

\begin{lema}\label{LemaCondicionAsterImplicaGPS}
If property $(\ast)$ is satisfied, then $\tilde \cF$ and $\tilde
\cF^\perp$ have a global product structure.
\end{lema}

\dem Consider any point $x_0 \in \tilde M$ and consider the set $G
= \{ z \in \tilde M \ : \ \tilde \cF^\perp (x_0) \cap \tilde
\cF(z) \neq \emptyset \}$. We have that $G$ is open from the local
product structure and by definition it is saturated by $\tilde \cF$. We must show that
$G$ is closed and since $\tilde M$ is connected this would
conclude.

Now, consider $z_0 \in \overline{G}$, using assumption ($\ast$) we
obtain that there exists $[y^-]< [z_0] <  [y^+]$ such that every
point $z$ such that $[y^-]<[z]<[y^+]$ verifies that its unstable
leaf intersects both $\tilde \cF(y^-)$ and $\tilde \cF(y^+)$.

Since $z_0\in \overline{G}$ we have that there are points $z_n \in
G$ such that $z_n \to \tilde \cF(z_0)$.

We get that eventually,  $[y^-] < [z_k] < [y^+]$ and thus we
obtain that there is a point $y \in \tilde \cF^\perp(x_0)$
verifying that $[y^-]<[y]<[y^+]$. We get that every leaf between
$\tilde \cF(y^-)$ and $\tilde \cF(y^+)$ is contained in $G$ from
assumption ($\ast$). In particular, $z_0\in G$ as desired. \lqqd

We must now show that property ($\ast$) is verified. To this end,
we will need the following lemma, which will allow us to show that
in the universal cover, the holonomy between a leaf of $\tilde
\cF$ and a translate by a deck transformation is bounded. This is
essential to obtain the proof and it is the place where the extent
of being "almost without holonomy" (which is measured by the value
of $K$) is used.

\begin{lema}\label{LemaConLongitudLCortoTodas}
Under the assumptions of Theorem
\ref{Teorema-EstructuraProductoGlobalMIO}, for sufficiently large
$K>0$ (which depends only on $M$ and $\delta$), there exists
$\ell>0$ such that every segment of $\cF^\perp(x)$ of length
$\ell$ intersects every leaf of $\cF$.
\end{lema}

We postpone the proof of this lemma to the next subsection
\ref{SubSectionPruebaLema}.

\demo{ of Theorem \ref{Teorema-EstructuraProductoGlobalMIO}} We
must prove that condition ($\ast$) is verified. Consider an
arbitrary point $z_0 \in \tilde M$.

We consider $\delta$ given by the size of local product structure
boxes and by Lemma \ref{LemaConLongitudLCortoTodas} we get a value
of $\ell>0$ such that every segment of $\cF^\perp$ of length
$\ell$ intersects every leaf of $\cF$.

As a consequence, by the pigeonhole principle, there exists $k>0$
such that every curve inside $\cF^\perp$ of length $k\ell$
verifies that it has a subarc whose endpoints are $\delta$-close
and joined by a curve in $\cF^\perp$ of length larger than $\ell$
(so, intersecting every leaf of $\cF$).

Since deck transformations are in one to one correspondence with
free homotopy classes of loops, there are finitely many deck
transformations which are represented by loops of length smaller
than $k\ell + \delta$. Here we are using the fact that $\pi_1(M)$
is abelian. Let us call $\Gamma_+^\trans$ to the (finite set) of
loops that are represented by transverse and positively oriented
loops of length smaller than or equal to $k \ell + \delta$ and
larger than $\ell$.

We obtain that in the one hand $\Gamma_+^\trans$ is non-empty and
finite from the argument above. On the other hand, because of
Lemma 6.4, we know that no element of $\Gamma_+^\trans$ fixes an
element of $\cL$ the space of leafs. Indeed, each element of
$\Gamma_+^\trans$ intersects every leaf of $\cF$ which implies
that its associated deck transformation must translate the lift of
every leaf in the universal cover in the same direction.

We deduce that there exists $\gamma_0 \in \Gamma_+^\trans$ such
that for every $\gamma \in \Gamma_+^\trans$ we have

$$  [z_0] < [z_0+\gamma_0] \leq [z_0+\gamma]$$

Consider now an arbitrary point $z \in \tilde \cF(z_0)$.

Let $\tilde \eta_z$ be the segment in $\tilde \cF^\perp_+(z)$ of
length $k \ell$ with one extreme in $z$. We can project $\tilde
\eta_z$ to $M$ and we obtain a segment $\eta_z$ transverse to
$\cF$ which contains two points $z_1$ and $z_2$ at distance
smaller than $\delta$ and such that the segment from $z_1$ to
$z_2$ in $\eta_z$ intersects every leaf of $\cF$. We denote
$\tilde z_1$ and $\tilde z_2$ to the lift of those points to
$\tilde \eta_z$.

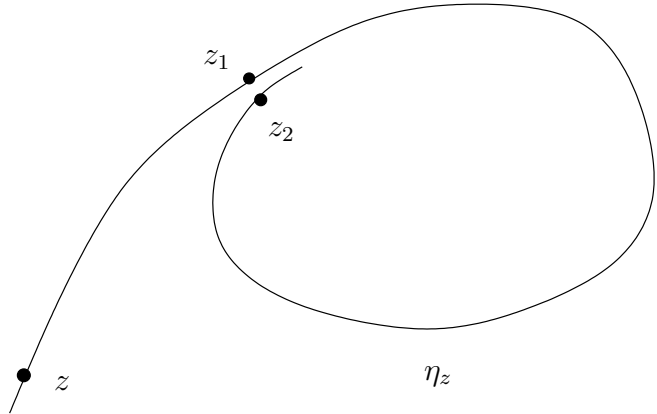
\begin{figure}[ht]
\begin{center}
\input{curvaeta.pstex_t}
\caption{\small{The curve $\eta_z$.}} \label{FiguraCurvaEta}
\end{center}
\end{figure}

Using the local product structure, we can modify slightly $\eta_z$
in order to create a closed curve $\eta_z'$ through $z_1$ which is
contained in $\eta_z$ outside $B_\delta(z_1)$, intersects every
leaf of $\cF$ and has length smaller than $k\ell + \delta$. The
lift of this curve clearly remains transverse to $\tilde \cF$ in
the positive direction.

Since $\eta_z$ essentially contains a loop of length smaller than
$k\ell + \delta$ we have that $\tilde \eta_z$ connects $[z_0]$
with $[\tilde z_1 + \gamma]$ where $\gamma$ belongs to
$\Gamma_+^\trans$. Moreover, since from $z$ to $\tilde z_1$ there
is a positively oriented arc of $\tilde \cF^\perp$ we get that
$[z_0]=[z] \leq [\tilde z_1]$ (notice that it is possible that
$z=\tilde z_1$).

From the choice of $\gamma_0$ above we deduce

$$[z_0] < [z_0 + \gamma_0] \leq [z_0 +\gamma] \leq [\tilde z_1 + \gamma] $$

\noindent where the last inequality follows from the fact that
$[z_0] \leq [\tilde z_1]$ and that the order is invariant by deck
transformations.

In particular, this proves that $\tilde \eta_z$ intersects the
leaf of $z_0 + \gamma_0$.

We have obtained that for $y^+= z_0 + \gamma_0$ there exists $L=
k\ell>0$ such that for every point $z\in \tilde \cF(z_0)$  the
segment of $\tilde \cF^\perp_+(z)$ of length $L$ intersects
$\tilde \cF(y^+)$.

This defines a continuous injective map from $\tilde \cF(z_0)$ to
$\tilde \cF(y^+)$ (injectivity follows from Proposition \ref{Proposicion-ArgumentoHaefliger}). Since the length of
the curves defining the map is uniformly bounded, this map is
proper and thus, a homeomorphism. The same argument applies to any
leaf $\tilde \cF(z_1)$ such that $[z_0] \leq [z_1] \leq [y^+]$.

For any $z_1$ such that $[z_0]\leq [z_1] \leq [y^+]$ we get that
$\tilde \cF^\perp(z_1)$ intersects $\tilde \cF(z_0)$. Since the
map defined above is a homeomorphism, we get that also $\tilde
\cF^\perp(z_0) \cap \tilde \cF(z_1) \neq \emptyset$.

A symmetric argument allows us to find $y^-$ with similar
characteristics. Using the fact that intersecting with leaves of
$\tilde \cF^\perp$ is a homeomorphism between any pair of leafs of
$\tilde \cF$ between $[y^-]$ and $[y^+]$ we obtain ($\ast$) as
desired.

Lemma \ref{LemaCondicionAsterImplicaGPS} finishes the proof. \lqqd

\subsection{Proof of Lemma \ref{LemaConLongitudLCortoTodas}}\label{SubSectionPruebaLema}

We first prove the following Lemma which allows us to bound the
topology of $M$ in terms of coverings of size $\delta$. Notice
that we are implicitly using that $\pi_1(M)$ as before to be able
to define a correspondence between (free) homotopy classes of
loops with elements of $\pi_1(M)$.

\begin{lema}\label{LemmaCotadePiunoenfunciondedeltacubrimiento}
Given a covering $\{V_1, \ldots , V_k\}$ of $M$ by contractible
open subsets there exists $K>0$ such that if $\eta$ is a loop in
$M$ intersecting each $V_i$ at most once\footnote{More precisely,
if $\eta$ is $\eta:[0,1] \to M$ with $\eta(0)=\eta(1)$ this means
that $\eta^{-1}(V_i)$ is connected for every $i$.}, then $[\eta]
\in \pi_1(M)$ has norm less than $K$.
\end{lema}

\dem We can consider the lift $p^{-1}(V_i)$ to the universal cover
of each $V_i$ and we have that each connected component of
$p^{-1}(V_i)$ has bounded diameter since they are simply connected
in $M$. Let $C_V> 0$ be a uniform bound on those diameters.

Let $K$ be such that every loop of length smaller than $2 k C_V$
has norm less than $K$ in $\pi_1(M)$.

Now, consider a loop $\eta$ which intersects each of the open sets
$V_i$ at most once. Consider $\eta$ as a function $\eta: [0,1] \to
M$ such that $\eta(0)=\eta(1)$. Consider a lift $\tilde \eta:
[0,1] \to M$ such that $p(\tilde \eta(t)) = \eta(t)$ for every
$t$.

We claim that the diameter of the image of $\tilde \eta$ cannot
exceed $k C_V$. Otherwise, this would imply that $\eta$ intersects
some $V_i$ more than once. Now, we can homotope $\tilde \eta$
fixing the extremes in order to have length smaller than $2 kC_V$.
This implies the Lemma.

\lqqd

Given $\delta$ of the uniform local product structure, we say that
a loop $\eta$ is a $\delta$-\emph{loop} if it is transverse to
$\cF$ and consists of a segment of a leaf of $\cF^\perp$ together
with a curve of length smaller than $\delta$.

\begin{lema}\label{LemaAbiertosSaturadosDeF}
There exists $K\geq 0$ such that if $O \en M$ is an open
$\cF$-saturated set such that $O \neq M$. Then, there is no
$\delta$-loop contained in $O$.
\end{lema}

\dem For every point $x$ consider $N_x = B_\delta(x)$ with
$\delta$ the size of the local product structure boxes. We can
consider a finite subcover $\{N_{x_1}, \ldots, N_{x_n} \}$ for
which Lemma \ref{LemmaCotadePiunoenfunciondedeltacubrimiento}
applies giving $K>0$.

Consider, an open set $O \neq M$ which is $\cF$-saturated. We must
prove that $O$ cannot contain a $\delta$-loop.

Let $\tilde O_0$ a connected component of the lift $\tilde O$ of
$O$ to the universal cover $\tilde M$. We have that the boundary
of $\tilde O_0$ consists of leaves of $\tilde \cF$ and if a
translation $\gamma \in \pi_1(M)$ verifies that

$$ \tilde O_0 \cap \gamma \tilde O_0 \neq \emptyset $$

\noindent then we must have that $\tilde O_0 = \gamma + \tilde
O_0$. This implies that $\gamma$ fixes the boundary leafs of
$\tilde O_0$: This is because the leaf space $\cL= \tilde
M/_{\tilde \cF}$ is homeomorphic to $\RR$ so that $\tilde O_0$
being connected and $\tilde \cF$ saturated is an open interval of
$\cL$. Since deck transformations preserve orientation, if they
fix an open interval then they must fix the boundaries.

The definition of $K$ then guarantees that if an element $\gamma$
of $\pi_1(M)$ makes $\tilde O_0$ intersect with itself, then
$\gamma$ must be larger than $K$. In particular, any $\delta$-loop
contained in $O$ must represent an element of $\pi_1(M)$ of length
larger than $K$.

Now consider a $\delta$-loop $\eta$. Proposition
\ref{Proposicion-ArgumentoHaefliger} implies that $\eta$ is in the
hypothesis of Lemma
\ref{LemmaCotadePiunoenfunciondedeltacubrimiento}. We deduce that
$\eta$ cannot be entirely contained in $O$ since otherwise its
lift would be contained in $\tilde O_0$ giving a deck
transformation $\gamma$ of norm less than $K$ fixing $\tilde O_0$
a contradiction.



\lqqd

\begin{cor}\label{CorolarioLoopSonTotales}
For the $K\geq 0$ obtained in the previous Lemma, if $\eta$ is a
$\delta$-loop then it intersects every leaf of $\cF$.
\end{cor}

\dem The saturation by $\cF$ of $\eta$ is an open set which is
$\cF$-saturated by definition. Lemma
\ref{LemaAbiertosSaturadosDeF} implies that it must be the whole
$M$ and this implies that every leaf of $\cF$ intersects $\eta$.
\lqqd

\demo{ of Lemma \ref{LemaConLongitudLCortoTodas}}  Choose $K$ as
in Lemma \ref{LemaAbiertosSaturadosDeF}. Considering a covering
$\{V_1, \ldots, V_k\}$ of $M$ by neighborhoods with local product
structure between $\cF$ and $\cF^\perp$ and of diameter less than
$\delta$.

There exists $\ell_0>0$ such that every oriented unstable curve of
length larger than $\ell_0$ traverses at least one of the $V_i's$.
Choose $\ell> (k+1)\ell_0$ and we get that every curve of length
larger than $\ell$ must intersect some $V_i$ twice in points say
$x_1$ and $x_2$. By changing the curve only in $V_i$ we obtain a
$\delta$-loop which will intersect the same leafs as the initial
arc joining $x_1$ and $x_2$.

Corollary \ref{CorolarioLoopSonTotales} implies that the mentioned
arc must intersect all leafs of $\cF$.

\lqqd

\subsection{Consequences of the global product structure in $\TT^3$}\label{SubSection-ConsequenciasGPS}

\subsubsection{} We say that a foliation $\cF$ in a Riemannian
manifold $M$ is \emph{quasi-isometric} if there exists $a,b \in
\RR$ such that for every $x,y$ in the same leaf of $\cF$ we have:

$$ d_{\cF}(x,y) \leq a d(x,y) + b $$

\noindent where $d$ denotes the distance in $M$ induced by the
Riemannian metric and $d_\cF$ the distance induced in the leaves
of $\cF$ by restricting the metric of $M$ to the leaves of $\cF$.

\begin{prop}\label{PropositionQuasiIsometria}
Let $\cF$ be a codimension one foliation of $\TT^3$ and
$\cF^\perp$ a transverse foliation. Assume the foliations $\tilde
\cF$ and $\tilde \cF^\perp$ lifted to the universal cover have
global product structure. Then, the foliation $\tilde \cF^\perp$
is quasi-isometric. Moreover, if $P$ is the plane given by Theorem
\ref{PropObienHayTorosObienTodoEsLindo}, there exists a cone $\cE$
transverse to $P$ in $\RR^3$ and $K>0$ such that for every $x\in
\RR^3$ and $y\in \tilde \cF^\perp(x)$ at distance larger than $K$
from $x$ we have that $y-x$ is contained in the cone $\cE$.
\end{prop}

\dem Notice that the global product structure implies that $\cF$
is Reebless. Let $P$ be the plane given by Theorem
\ref{PropObienHayTorosObienTodoEsLindo}.

Consider $v$ a unit vector perpendicular to $P$ in $\RR^3$.

\begin{af} For every $N>0$ there
exists $L$ such that every segment of $\tilde \cF^\perp_+$ of
length $L$ centered at a point $x$ intersects both $P+x+Nv$ and
$P+x-Nv$.
\end{af}

\dem If this was not the case, we could find arbitrarily large
segments $\gamma_k$ of leaves of $\tilde \cF^\perp$ centered at a
point $x_x$ with length larger than $k$ and such that they do not
intersect either $P+x_k +Nv$ or $P+x_k -Nv$. Without loss of
generality, and by taking a subsequence we can assume that they do
not intersect $P+x_k +Nv$.

Since the foliations are invariant under translations, we can
assume that the sequence $x_k$ is bounded and by further
considering a subsequence, that $x_k \to x$.

We deduce that $\tilde \cF^{\perp}(x)$ cannot intersect $P+x +
(N+1)v$ which in turn implies (by Theorem
\ref{PropObienHayTorosObienTodoEsLindo}) that $\tilde
\cF^\perp(x)$ cannot intersect the leaf of $\tilde \cF$ through
the point $x+ (N+1 +2R)v$ contradicting global product structure.

\finobs

This implies quasi-isometry since having length larger than $kL$
implies that the endpoints are at distance at least $kN$.

It also implies the second statement since assuming that it does
not hold, we get a sequence of points $x_n$, $y_n$ such that the
distance is larger than $n$ and such that the angle between $y_n -
x_n$ with $P$ is smaller than $1/n$. This implies that the length
of the arc of $\tilde \cF^\perp$ joining $x_n$ and $y_n$ is larger
than $n$ and that it does not intersect $P + x_n + (n
\sin(\frac{1}{n})+2R) v$ contradicting the claim.

\lqqd


\section{Partially hyperbolic diffeomorphisms isotopic to Anosov}\label{Section-PHAnosov}

\subsection{Preliminaries and notation} \subsubsection{} In this section we give a proof of Theorem A.

We shall assume that $f: \TT^3 \to \TT^3$ is an almost dynamical
coherent partially hyperbolic diffeomorphism isotopic to a linear
Anosov automorphism $A: \TT^3 \to \TT^3$ with splitting of the
form $T \TT^3 = E^{cs} \oplus E^u$ with $\dim E^u=1$.

This means that $f_\ast$ coincides with the lift of $A$ to $\RR^3$
and is a hyperbolic matrix. By abuse of notation $A$ will
denote both the hyperbolic matrix in $\RR^3$ and the
diffeomorphism of $\TT^3$.  We will denote as $\cF$ to the foliation given by the
definition of almost dynamical coherent which we know is Reebless.

\subsubsection{} It is important to remark that we are not assuming that the stable
dimension of $A=f_\ast$ coincides with the one of $E^{cs}$. In
fact, many of the arguments below become much easier in the case
$A$ has stable dimension $2$. 

\subsection{A planar direction for the foliation transverse to $E^u$}

\subsubsection{} Proposition \ref{PropExisteSemiconjugacion} implies the existence
of a continuous $\ZZ^3$-periodic surjective function $H: \RR^3\to
\RR^3$ which verifies
$$H \circ \tilde f = A \circ H$$

\noindent and such that $d(H(x),x)< K_1$ for every $x\in \RR^3$.

We can prove:

\begin{lema}\label{LemaHdemediaFuesNoacotado} For every $x\in \RR^3$ we have that $H(\tilde \cF^u_+(x))$ is unbounded.
\end{lema}
\dem Otherwise, for some $x\in \RR^3$, the unstable leaf $\tilde
\cF^u_+(x)$ would be bounded. Since its length is infinite one can
find two points in $\tilde \cF^u_+(x)$ in different local unstable
leafs at arbitrarily small distance. This contradicts Corollary
\ref{CorolarioConsecuenciasReeb} (i). \lqqd

\subsubsection{} Since $\cF$ is transverse to the unstable direction, we get by
Theorem \ref{CorolarioConsecuenciasReeb} that it is a Reebless
foliation so that we can apply Theorem
\ref{PropObienHayTorosObienTodoEsLindo}. We intend to prove in
this section that option (ii) of this Theorem
\ref{PropObienHayTorosObienTodoEsLindo} is not possible when $f$
is isotopic to Anosov (see the example in \cite{HHU} where that
possibility occurs).

Notice that if we apply $\tilde f^{-1}$ to the foliation $\tilde
\cF$, then the new foliation $f^{-1}(\tilde \cF)$ is still
transverse to $E^u$ so that Theorem
\ref{PropObienHayTorosObienTodoEsLindo} still applies. In fact, if
$P$ is the plane obtained for $\tilde \cF$, then the plane which
the proposition will give for $\tilde f^{-1}(\tilde \cF)$ will be
$A^{-1}(P)$: This is immediate by recalling that $\tilde f$ and
$A$ are at bounded distance while two planes which are not
parallel have points at arbitrarily large distance.

\subsubsection{} 
We say that a subspace $P$ is \emph{almost parallel} to a
foliation $\tilde \cF$ if there exists $R>0$ such that for every
$x\in \RR^3$ we have that $P+x$ lies in an $R$-neighborhood of
$\tilde \cF(x)$ and $\tilde \cF(x)$ lies in a $R$-neighborhood of
$P+x$.

\begin{prop}\label{PropSifIsotopicoAAnosovNoHayTorosEnFoliacion}
Let $f: \TT^3 \to \TT^3$ be a partially hyperbolic diffeomorphism
of the form $T\TT^3 = E^{cs} \oplus E^u$ (with $\dim E^{cs}=2$)
isotopic to a linear Anosov automorphism and $\cF$ a foliation
transverse to $E^u$. Then, there exists a two dimensional subspace
$P \en \RR^3$ which is almost parallel to $\tilde \cF$.
\end{prop}


\dem Assume by contradiction that option (ii) of Proposition
\ref{PropObienHayTorosObienTodoEsLindo} holds. Then, there exists
a plane $P \en \RR^3$ whose projection to $\TT^3$ is a two
dimensional torus and such that every leaf of $\tilde \cF^u$,
being transverse to $\tilde \cF$, remains at bounded
distance from $P$.

Since $f$ is isotopic to a linear Anosov automorphism $A$ we know
that $P$ cannot be invariant under $A$ (see Proposition
\ref{PropAnosovSonIrreducibles}). So, we have that the intersection between $P$, $A^{-1}(P)$ and $A^{-2}(P)$ is a unique point (the origin).



We get that for every point $x$ we have that $\tilde \cF^u_+(x)$
must lie within bounded distance from $P$ as well as from
$A^{-1}(P)$ (since when we apply $\tilde f^{-1}$ to $\tilde \cF$
the leaf close to $P$ becomes close to $A^{-1}(P)$) and $A^{-2}(P)$. This implies that in
fact $\tilde \cF^u_+(x)$ lies within bounded distance from the origin, but this contradicts the fact that $\tilde \cF^u_+(x)$ is unbounded (Lemma \ref{LemaHdemediaFuesNoacotado}). \lqqd

\subsection{Global Product Structure}\label{SectionGPS}

%

\subsubsection{}
Proposition \ref{PropSifIsotopicoAAnosovNoHayTorosEnFoliacion}
implies that the foliation $\tilde \cF$ is quite well behaved. In
this section we shall show that the properties we have showed for
the foliations and the fact that $\tilde \cF^u$ is $\tilde
f$-invariant while the foliation $\tilde \cF$ remains with a
uniform local product structure with $\tilde \cF^u$ when iterated
backwards (see Lemma \ref{LemaEstructuraProductoLocalUniforme})
implies that there is a global product structure. 

The main result of this section is thus the following:

\begin{prop}\label{PropositionGlobalProductStructure}
Given $x, y\in \RR^3$ we have that $\tilde \cF(x) \cap \tilde
\cF^u(y) \neq \emptyset$. This intersection consists of exactly
one point.
\end{prop}

Notice that uniqueness of the intersection point follows directly
from Theorem \ref{CorolarioConsecuenciasReeb} (i). The proof
consists in showing that for sufficiently large $n$ we have that
$f^{-n}(\cF)$ and $\cF^u$ are in the conditions of Theorem
\ref{Teorema-EstructuraProductoGlobalMIO}.

\subsubsection{} 
We start by proving a result which gives that the size of local
product structure boxes between $f^{-n}(\cF)$ and $\cF^u$ can be
chosen independent of $n$. We shall denote as $\DD^2= \{ z\in \CC
\ : \ |z|\leq 1 \}$.

\begin{lema}\label{LemaEstructuraProductoLocalUniforme}
There exists $\delta>0$ such that for every $x\in \RR^3$ and
$n\geq 0$ there exists a closed neighborhood $V_x^n$ containing
$B_\delta(x)$ such that it admits $C^0$-coordinates $\varphi_x^n:
\DD^2 \times [-1,1] \to \RR^3$ such that:
\begin{itemize}
\item[-] $\varphi_x^n(\DD^2\times [-1,1]) = V_x^n$ and
$\varphi_x^n(0,0)= x$. \item[-] $\varphi_x^n(\DD^2 \times \{t\}) =
\tilde f^{-n}(\tilde \cF(\tilde f^n(\varphi_x^n(0,t)))) \cap
V_x^n$ for every $t \in [-1,1].$ \item[-] $ \varphi_x^n(\{s\}
\times [-1,1]) = \tilde \cF^u(\varphi_x^n(s,0)) \cap V_x^n$ for
every $s\in \DD^2.$
\end{itemize}
\end{lema}


\dem  Notice first that  $f^{-n}(\cF)$ is tangent to a cone
transverse to $E^u$ and independent of $n$. Let us call this cone
$\cE^{cs}$.

Given $\epsilon>0$ we can choose a neighborhood $V_{\epsilon}$ of
$x$ contained in $B_{\epsilon}(x)$ such that the following is
verified:

\begin{itemize}
\item[-] There exists a two dimensional disk $D$ containing $x$
such that $V_{\epsilon}$ is the union of segments of $\cF^u(x)$ of
length $2\epsilon$ centered at points in $D$. This defines two
boundary disks $D^+$ and $D^-$ contained in the boundary of
$V_\eps$. \item[-] By choosing $D$ small enough, we get that there
exists $\epsilon_1>0$ such that every curve of length $\epsilon_1$
starting at a point $y \in B_{\epsilon_1}(x)$ tangent to
$\cE^{cs}$ must leave $V_{\epsilon}$ and intersects $\partial
V_{\epsilon}$ in $\partial V_{\epsilon} \setminus (D^+ \cup D^-)$.
\end{itemize}

Notice that both $\epsilon$ and $\epsilon_1$ can be chosen
uniformly in $\RR^3$ because of compactness of $\TT^3$ and uniform
transversality of the foliations.

This implies that every disk of radius $\epsilon$ tangent to
$\cE^{cs}$ centered at a point $z \in B_{\epsilon_1}(x)$ must
intersect the unstable leaf of every point in $D$, in particular,
there is a local product structure of uniform size around each
point in $\RR^3$.

Now, we can choose a continuous chart (recall that the foliations
are only continuous, but with $C^1$-leaves) around each point
which sends horizontal disks into disks transverse to $E^u$ and
vertical lines into leaves of $\tilde \cF^u$ containing a fixed
ball around each point $x$ independent of $n\geq 0$ giving the
desired statement. \lqqd

\begin{obs}\label{RemarkEPLUniforme}
We obtain that there exists $\eps>0$ such that for every $x\in
\RR^3$ there exists $V_x \en \bigcap_{n\geq 0} V_x^n$ containing
$B_\eps(x)$ admitting $C^1$-coordinates $\psi_x : \DD^2 \times
[-1,1] \to \RR^3$ such that:
\begin{itemize}
\item[-] $\psi_x (\DD^2 \times [-1,1]) = V_x$ and $\psi_x(0,0)=x$.
\item[-] If we consider $V_x^\eps= \psi_x^{-1}(B_\eps(x))$ then
one has that for every $y \in V_x^\eps$ and $n\geq 0$ we have
that: $$ \psi_x^{-1}(\tilde f^{-n}(\tilde \cF(\tilde f^n(y))) \cap
V_x) $$
    \noindent is the graph of a function $h_y^n: \DD^2 \to [-1,1]$ which has uniformly bounded derivative in $y$ and $n$.
\end{itemize}
Indeed, this is given by considering a $C^1$-chart $\psi_x$ around
every point such that its image covers the $\eps$-neighborhood of
$x$ and sends the $E$-direction to an almost horizontal direction
and the $E^u$-direction to an almost vertical direction. See for
example \cite{BuW} section 3 for more details on this kind of
constructions. \finobs
\end{obs}

\subsubsection{} The next lemma shows that after iterating the foliation backwards, one
gets that it becomes nearly irrational so that we can apply
Theorem \ref{Teorema-EstructuraProductoGlobalMIO}.

\begin{lema}\label{LemaIterandoParaAtrasSeHaceIrracional}
Given $K>0$ there exists $n_0>0$ such that for every $x\in \RR^3$
and for every $\gamma \in \ZZ^3$ with norm less than $K$ we have
that $$\tilde f^{-n_0}(\tilde \cF(x))+ \gamma \neq \tilde
f^{-n_0}(\tilde \cF(x)) \qquad  \forall x \in \RR^3.$$
\end{lema}

\dem Notice that $\tilde f^{-n}(\tilde \cF)$ is almost parallel to
$A^{-n}(P)$. Notice that $A^{-n}(P)$ has a converging subsequence
towards a totally irrational plane $\tilde P$ (see Remark
\ref{RemarkSubespaciosInvariantesAnosov}).

We can choose $n_0$ large enough such that no element of $\ZZ^3$
of norm smaller than $K$ fixes $A^{-n_0}(P)$.

Notice first that $\tilde f^{-n_0}(\tilde \cF)$ is almost parallel
to $A^{-n_0}(P)$. Now, assuming that there
is a translation $\gamma$ which fixes a leaf of $\tilde
f^{-n_0}(\tilde \cF(x))$ we get that the leaf $p(\tilde
f^{-n_0}(\tilde \cF(x)))$ contains a loop homotopic to $\gamma$.

This implies that it is at bounded distance from the line which is
the lift of the canonical (linear) representative of $\gamma$.
This implies that $\gamma$ fixes $A^{-n_0}(P)$ and thus has norm
larger than $K$ as desired (see also Lemma
\ref{RemarkHojasCerradas}).

\lqqd
%

\demo{of Proposition \ref{PropositionGlobalProductStructure}} By
Theorem \ref{CorolarioConsecuenciasReeb} we know that all the
leaves of $\tilde \cF$ are simply connected. Proposition
\ref{RemarkElCocienteDeLaFoliacionEsR} implies that the leaf space
of $\tilde \cF$ is homeomorphic to $\RR$. All this properties
remain true for the foliations $\tilde f^{-n}(\tilde \cF)$ since
they are diffeomorphisms at bounded distance from linear
transformations.

Lemma \ref{LemaEstructuraProductoLocalUniforme} gives that the
size of the local product structure between $\tilde f^{-n}(\tilde
\cF)$ and $\tilde \cF^u$ does not depend on $n$.

Using Lemma \ref{LemaIterandoParaAtrasSeHaceIrracional} we get
that for some sufficiently large $n$ the foliations $\tilde
f^{-n}(\tilde \cF)$ and $\tilde \cF^u$ are in the hypothesis of
Theorem \ref{Teorema-EstructuraProductoGlobalMIO} which gives
global product structure between $\tilde f^{-n}(\tilde \cF)$ and
$\tilde \cF^u$. Since $\tilde \cF^u$ is $\tilde f$-invariant and
$f$ is a diffeomorphism we get that there is a global product
structure between $\tilde \cF$ and $\tilde \cF^u$ as desired.

\lqqd


\begin{cor}\label{PropositionQuasiIsometriaII}
The foliation $\tilde \cF^u$ is quasi-isometric. Moreover, there
exist one dimensional subspaces $L_1$ and $L_2$ of $E^u_A$
transverse to $P$ and $K>0$ such that for every $x\in \RR^3$ and
$y\in \tilde \cF^u(x)$ at distance larger than $K$ from $x$ we
have that $H(y)-H(x)$ is contained in the cone of $E^u_A$ with
boundaries $L_1$ and $L_2$ and transverse to $P$.
\end{cor}

Notice that if $A$ has stable dimension $2$ then $L_1=L_2 =E^u_A$.

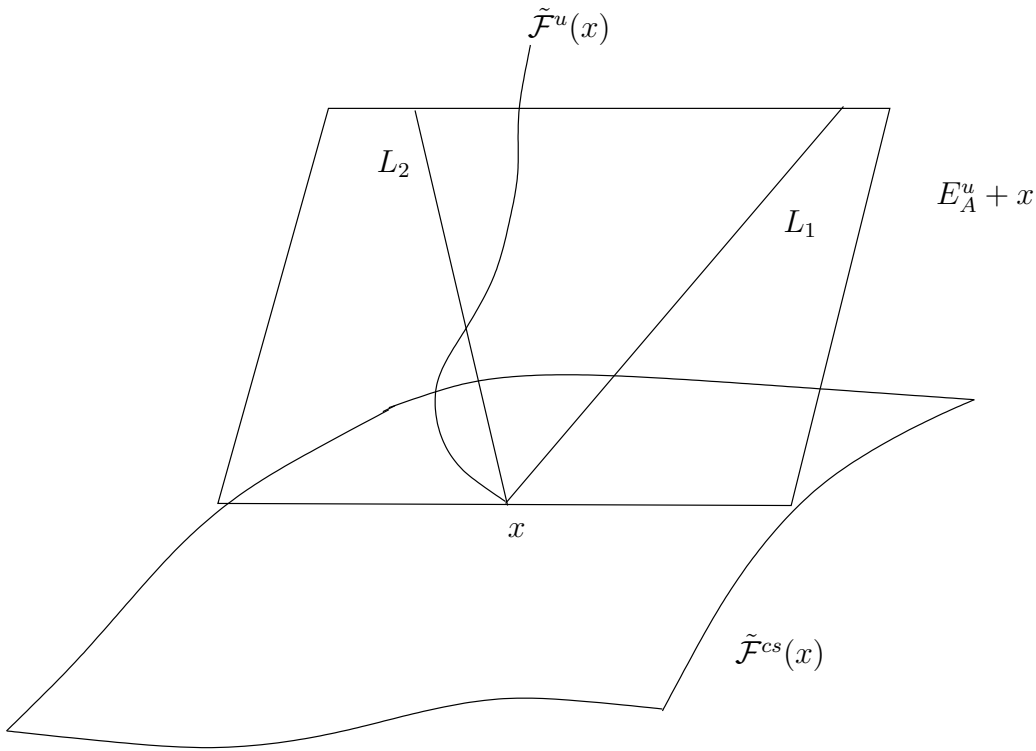
\begin{figure}[ht]
\begin{center}
\input{cono.pstex_t}
\caption{\small{The unstable leaf of $x$ remains close to the cone
bounded by $L_1$ and $L_2$.}} \label{FiguraCono}
\end{center}
\end{figure}

\dem This is a direct consequence of Proposition
\ref{PropositionQuasiIsometria} and the fact that the image of
$\tilde \cF^u(x)$ by $H$ is contained in $E^u_A +H(x)$. \lqqd

\begin{obs}\label{RemarkHesInjectivaEnInestables}
Since points which are sent to the same point by $H$ must have
orbits remaining at bounded distance, the quasi-isometry of
$\tilde \cF^u$ implies that $H$ must be injective on leaves of
$\tilde \cF^u$. \finobs
\end{obs}

\subsection{Complex eigenvalues}\label{SubSectionComplejosYQIsom}

The following proposition has interest only in the case $A$ has
stable dimension one. It establishes the last statement of Theorem
A.

\begin{prop}\label{PropositionAnotieneinestablesComplejos} The matrix $A$ cannot have complex unstable eigenvalues.
\end{prop}

\dem Assume that $A$ has complex unstable eigenvalues, in
particular $E^u_A$ is two-dimensional. Consider a fixed point
$x_0$ of $\tilde f$.

Recall that by Lemma \ref{LemaHdemediaFuesNoacotado} the set $\eta
=H(\tilde \cF^u_+(x_0))$ is an unbounded continuous curve in
$E^u_A$. Since $x_0$ is fixed and since $H$ is a semiconjugacy, we
have that $\eta$ is $A$-invariant.

On the other hand, by Corollary \ref{PropositionQuasiIsometriaII}
we have that $\eta$ is eventually contained in a cone between two
lines $L_1$ and $L_2$.

This implies that $A$ cannot have complex unstable eigenvalues
(recall that they should have irrational angle by Proposition
\ref{PropAnosovSonIrreducibles}) since a matrix which preserves an
unbounded connected subset of a cone cannot have complex
eigenvalues with irrational angle. \lqqd

\subsection{Dynamical Coherence}\label{SubSectionDynamicalCoherence}

\subsubsection{} In this section we shall show dynamical coherence of almost
dynamically coherent partially hyperbolic diffeomorphisms isotopic
to linear Anosov automorphisms (and uniqueness of the
$f$-invariant foliation tangent to $E^{cs}$). This will complete
the proof of Theorem A.

The proof of dynamical coherence is quite lengthy and technical,
so let us first explain the strategy and the potential
difficulties we must deal with in the proof.

The idea is quite direct, consider the foliation $\cF$ transverse
to $E^u$ and iterate it backwards hoping that in the limit it will
converge to the desired $f$-invariant foliation $\cW^{cs}$ tangent
to $E^{cs}$. However, the following difficulties might appear:

\begin{itemize}
\item The leaves of $f^{-n}(\cF)$ become uniformly tangent to
$E^{cs}$ (this is a standard graph transform argument), however,
it is not clear that their limit is unique nor that in the limit
different leafs might not merge. To overcome this problem, we use
the semiconjugacy to show that each limit disk by the backward
iteration is sent by the semiconjugacy to a certain $A$-invariant
plane. This is shown in the first claim inside the proof.

\item The second difficulty is that it is not clear that this
invariant plane remains transverse to the cone where we know the
strong unstable manifolds are, so that in the limit we could in
principle loose the global product structure (see Question
\ref{QuestionPlane} below). Moreover, it must be shown that the
limit of each leaf is sent to a different translate of the
$A$-invariant plane. This is shown by using the second claim in
the last claim of the proof. By working with local plaques we can
treat both difficulties at once (the fact that one leaf might have
more than one limit and that many leafs might have limits which
merge).

\end{itemize}

\begin{teo}\label{TeoremaCoherencia}
Let $f:\TT^3 \to \TT^3$ be an almost dynamically coherent
partially hyperbolic diffeomorphism of the form $T\TT^3 = E^{cs}
\oplus E^u$ isotopic to a linear Anosov automorphism. Then, there
exists an $f$-invariant foliation $\cF^{cs}$ tangent to $E^{cs}$.
If $\tilde \cF^{cs}$ denotes the lift to $\RR^3$ of this
foliation, then $H(\tilde \cF^{cs}(x)) = P^{cs} + H(x)$ where
$P^{cs}$ is an $A$-invariant subspace and $E^u_A$ is not contained
in $P^{cs}$.
\end{teo}

\dem Consider the foliation $\tilde \cF$, by Proposition
\ref{PropSifIsotopicoAAnosovNoHayTorosEnFoliacion} we have a plane
$P$ which is almost parallel to $\tilde \cF$.

Let $P^{cs}$ be the limit of $A^{-n}(P)$ which is an $A$-invariant
subspace. Since we have proved that $A$ has no complex unstable
eigenvalues (Proposition
\ref{PropositionAnotieneinestablesComplejos}) and since $P$ is
transverse to $E^u_A$ (Proposition
\ref{PropositionQuasiIsometria}), this plane is well defined (see
Remark \ref{RemarkSubespaciosInvariantesAnosov}).

Notice that the transversality of $P$ with $E^u_A$ implies that
$P^{cs}$ contains $E^s_A$, the eigenspace associated with stable
eigenvalues (in the case where $A$ has stable dimension $2$ we
thus have $P^{cs}=E^s_A$).

Since $P^{cs}$ is $A$-invariant, we get that it is totally
irrational so that no deck transformation fixes $P^{cs}$.

Using Remark \ref{RemarkEPLUniforme} we obtain $\eps>0$ such that
for every $x\in \RR^3$ there are neighborhoods $V_x$ containing
$B_\eps(x)$ admitting $C^1$-coordinates $\psi_x : \DD^2 \times
[-1,1] \to V_x$ such that:

\begin{itemize}
\item[-] For every $y\in B_\eps(x)$ we have that if we denote as
$W^x_n(y)$ to the connected component containing $y$ of $V_x \cap
\tilde f^{-n}(\tilde \cF(f^n(y)))$ then the set
$\psi_x^{-1}(W^x_n(y))$ is the graph of a $C^1$-function
$h^{x,y}_n: \DD^2 \to [-1,1]$ with bounded first
derivatives\footnote{Notice that the graphs are tangent to a
distribution contained in a fixed cone field which is transverse
to the direction of projection. See Remark
\ref{RemarkEPLUniforme}. This construction here can be compared
with classical graph transform constructions which do essentially
the same, see \cite{HPS} or \cite{BuW} section 3 where certain
\emph{fake foliations} are constructed in this way.}.
\end{itemize}

Using the fact that these graphs have bounded first derivative we
get that $\{h^{x,y}_n\}$ is pre-compact in the space of Lipschitz
functions from $\DD^2$ to $[-1,1]$. 

For every $y\in B_\eps(x)$ there exists $\cJ_y^x$ a set of indices
such that for every $\alpha \in \cJ_y^x$ we have a Lipschitz
function $h^{x,y}_{\infty,\alpha}: \DD^2 \to [-1,1]$ and $n_j \to
+\infty$ such that:

$$ h^{x,y}_{\infty,\alpha} = \lim_{j \to + \infty} h^{x,y}_{n_j} $$

The set of indices included in $\cJ_y^x$ is given by the fact
that, a priori, different subsequences might have different
limits\footnote{One can think of $\cJ_y^x$ as an equivalence
relation between the subsequences $n_j \to \infty$ such that the
functions $h^{x,y}_{n_j}$ converges and the equivalence is given
by having the same limit.}. The main point from now on is to show
that the limits are unique.

Every $h^{x,y}_{\infty,\alpha}$ gives rise to a graph whose image
by $\psi_x$ we denote as $W_{\infty,\alpha}^{x}(y)$. This
Lipschitz manifold verifies that it contains $y$ and is tangent to
a small cone\footnote{This is enough to show that the manifolds
are $C^1$ and tangent to $E^{cs}$ but we shall deduce this later
from invariance which is the standard way.} around $E^{cs}$.

\begin{af}
For every $z\in B_\eps(x)$ and every $\alpha \in \cJ_z^x$, we have
that $H(W^x_{\infty,\alpha}(z)) \en P^{cs}+H(z)$
\end{af}

\dem Consider $z \in B_\eps(x)$ and $\alpha \in \cJ_z^x$. One can
find $n_j \to \infty$ such that $W^x_{n_j}(z) \to
W^x_{\infty,\alpha}(z)$. Pick now any $y \in
W^x_{\infty,\alpha}(z)$.

In the coordinates $\psi_x$ of $V_x$, we can find a sequence
$z_{n_j} \in W^x_{n_j}(z) \cap \tilde \cF^u(y)$ such that $z_{n_j}
\to y$. Moreover, we have that $\tilde f^{n_j}(z_{n_j}) \in \tilde
\cF(\tilde f^{n_j}(z))$. Assume that $H(y) \neq H(z)$ (otherwise
there is nothing to prove).

We have, by continuity of $H$ that $H(z_{n_j})\to H(y) \neq H(z)$.
We must show that $H(y) - H(z) \in P^{cs}$.

Notice first that $\tilde f^{n_j}(z)$ and $\tilde f^{n_j}(z_j)$
lie in the same leaf of $\tilde \cF$. This implies that $\tilde
f^{n_j}(z) - \tilde f^{n_j}(z_{n_j})$, and consequently  $H(\tilde
f^{n_j}(z)) - H(\tilde f^{n_j}(z_j))$, are uniformly bounded in
the direction transverse to $P$, the plane defined above.

Using the semiconjugacy, this is the same as saying that
$A^{n_j}(H(z)) - A^{n_j}(H(z_{n_j}))$ is uniformly bounded in the
direction transverse to $P$.

Since $A^{-n}(P)$ converges exponentially to $P^{cs}$, and $A$ is
linear, this implies that given $R>0$ one has that
$A^{-n}(B_R(P))$ converges uniformly to $P^{cs}$. This in turn
implies that $H(y)-H(z)=\lim_j H(z_{n_j}) - H(z) \in P^{cs}$ as
desired.

\finobs

Assuming that $P^{cs}$ does not intersect the cone bounded by
$L_1$ and $L_2$ this finishes the proof since one sees that each
leaf of $\tilde \cF^u$ can intersect the pre-image by $H$ of
$P^{cs}+y$ in a unique point, thus showing that the partition of
$\RR^3$ by the pre-images of the translates of $P^{cs}$ defines a
$\tilde f$-invariant foliation (and also invariant under deck
transformations). We leave to the interested reader the task of
filling the details of the proof in this particular case, since we
will continue by giving a proof which works in all cases.

We will prove that $H$ cannot send unstable intervals into the
same plane parallel to $P^{cs}$.

\begin{af}
Given $\eta :[0,1] \to \tilde \cF^u(x)$ a non-trivial curve, we
have that $H(\eta([0,1]))$ is not contained in $P^{cs} +
H(\eta(0))$.
\end{af}

\dem We prove this claim by contradiction. Assume that
$H(\eta([0,1])) \en P^{cs} + H(\eta(0))$.

Consider $C_\eps$ given by Theorem
\ref{CorolarioConsecuenciasReeb} (iv) for $\eps$. Moreover,
consider $L$ large enough such that $C_\eps L > \Vol (\TT^3)$.

Since $\tilde \cF^u$ is $\tilde f$-invariant and $P^{cs}$ is
$A$-invariant we deduce that we can assume that the length of
$\eta$ is arbitrarily large, in particular larger than $2L$.

We will show that $H(B_\eps(\eta([a,b]))) \en P^{cs} + H(\eta(0))$
where $0< a < b < 1$ and the length of $\eta([a,b])$ is larger
than $L$.

Having volume larger than $\Vol(\TT^3)$ there must be a deck
transformation $\gamma \in \ZZ^3$ such that $\gamma +
B_\eps(\eta([a,b])) \cap B_\eps(\eta([a,b]))\neq \emptyset$. This
in turn gives that $\gamma +H(B_\eps(\eta([a,b]))) \cap
H(B_\eps(\eta([a,b]))) \neq \emptyset$ and thus $\gamma + P^{cs}
\cap P^{cs} \neq \emptyset$. Since $P^{cs}$ is totally irrational
this is a contradiction.

It remains to show that, under the assumption that $H(\eta([0,1]))
\en H(\eta(0)) + P^{cs}$ we also have that $H(B_\eps(\eta([a,b])))
\en P^{cs} + H(\eta(0))$. By the previous claim, we know that if
$z,w \in W^x_{\infty,\alpha}(y)$ for some $\alpha \in \cJ_y$, then
$H(z)-H(w) \in P^{cs}$.

Consider $a,b \in [0,1]$ such that $\tilde \cF^u(x) \cap
B_\eps(\eta([a,b])) \en \eta([0,1])$. By Theorem
\ref{CorolarioConsecuenciasReeb} we have that such $a,b$ exist and
we can choose them\footnote{To avoid confusions, we remark that
the reason why we must choose $a >0$ and $b< 1$ is that the
$\eps$-neigborhood of $\eta(0)$ intersects $\cF^u(x)$ in points
which are not in $\eta([0,1])$.} in order that the length of
$\eta([a,b])$ is larger than $L$.

Let $z\in B_\eps(\eta([a,b]))$ and choose $w \in \eta([a,b])$ such
that $z \in B_\eps(w)$. We get that for every $\alpha \in \cJ_z^w$
we have that $W^{w}_{\infty, \alpha}(z) \cap \eta([0,1]) \neq
\emptyset$. Since $H(\eta([0,1])) \en P^{cs} + H(\eta(0))$ and by
the previous claim, we deduce that $H(z) \en P^{cs} + H(\eta(0))$
finishing the proof.

\finobs

Now we are in conditions of showing that for every point $x$ and
for every point $y\in B_\eps(x)$ the limit $W^x_\infty(y)$ of the
manifolds $W^x_n(y)$ is unique and tangent to $E^{cs}$.


Indeed, assume first that the manifolds $W^x_n(y)$ have a unique
limit for every $x \in \RR^3$ and $y\in B_\eps(x)$ and that for
any pair points $y,z\in B_\eps(x)$ these limits are either
disjoint or equal (see the claim below). One has that the set of
manifolds $W^x_\infty(y)$ forms an $f$-invariant plaque family in
the following sense:

\begin{itemize}
\item[-] $\tilde f(W^x_\infty(y)) \cap W^{\tilde f(x)}_\infty
(\tilde f(y))$ is relatively open whenever $\tilde f(y) \in
B_\eps(\tilde f(x))$.
\end{itemize}

We must thus show that these plaque families form a foliation.
Consider $z,w \in B_\eps(x)$ we have that $W^x_\infty (z) \cap
\tilde \cF^u(w) \neq \emptyset$ and in fact consists of a unique
point (see Theorem \ref{CorolarioConsecuenciasReeb} (i)). Since
the intersection point varies continuously and using that plaques
are either disjoint or equal we obtain a continuous map from
$\DD^2 \times [-1,1]$ to a neighborhood of $x$ sending horizontal
disks into plaques. This implies that the plaques form an
$f$-invariant foliation as desired. The fact that this foliation
is tangent to $E^{cs}$ follows directly from the fact that its
leafs are tangent to a small cone around $E^{cs}$ and the
foliation is invariant (compare for example with the $C^r$-section
theorem 3.5 in \cite{HPS}, though this case is of course much
simpler).

It thus remains to show the following:

\begin{af}
Given $x \in \RR^3$ and $y,z \in B_\eps(x)$ we have that there is
a unique limit of $W^x_\infty (y)$ and $W^x_\infty(z)$ and they
are either disjoint or coincide. More precisely, for every $\alpha
\in \cJ^x_y$ and $\beta \in \cJ^x_z$ ($z$ could coincide with $y$)
we have that either $h^{x,y}_{\infty,\alpha}=
h^{x,z}_{\infty,\beta}$ or the graphs are disjoint.
\end{af}

\dem Assuming the claim does not hold, one obtains $y,z \in
B_\eps(x)$ such that $h^{x,y}_{\infty,\alpha}$ and
$h^{x,z}_{\infty,\beta}$ coincide at some point but whose graphs
are different for some $\alpha \in \cJ_y^x$ and $\beta \in
\cJ_z^x$. In particular, there exists a point $t \in \DD^2$ which
is in the boundary of where both functions coincide. We assume for
simplicity\footnote{If it were not the case we would need to
change the coordinates and perform the same proof, but to avoid
charging the notation we choose to make this (unnecessary)
assumption.} that $\psi_x (t)$ belongs to $B_\eps(x)$.

Let $\gamma : [0,1] \to B_\eps(x)$ be a non-trivial arc of $\tilde
\cF^u$ joining the graphs of $h^{x,y}_{\infty,\alpha}$ and
$h^{x,z}_{\infty,\beta}$. Since the graphs of both
$h^{x,y}_{\infty,\alpha}$ and $h^{x,z}_{\infty,\beta}$ separate
$V_x$ we have that every point $w \in \gamma((0,1))$ verifies that
for every $\delta \in \cJ_w^x$ one has that $W^x_{\infty,
\delta}(w)$ intersects at least one of $W^x_{\infty,\alpha}(y)$ or
$W^x_{\infty, \beta}(z)$. By the first claim we get that $H(w) \in
P^{cs} + H(y) = P^{cs}+ H(z)$ a contradiction with the second
claim.

\finobs

\lqqd

\subsubsection{} We can in fact obtain a stronger property since
our results allows us to show that in fact $E^{cs}$ is uniquely
integrable into a $f$-invariant foliation. There are stronger
forms of unique integrability (see \cite{BuW2} and \cite{BF}).

\begin{prop}\label{Prop-UnicidadDeFoliacionTgAEcs}
There is a unique $f$-invariant foliation $\cF^{cs}$  tangent to
$E^{cs}$. Moreover, the plane $P^{cs}$ given by Theorem
\ref{PropObienHayTorosObienTodoEsLindo} for this foliation is
$A$-invariant and contains the stable eigenspace of $A$.
\end{prop}

\dem Assume there are two different $f$-invariant foliations
$\cF^{cs}_1$ and $\cF^{cs}_2$ tangent to $E^{cs}$.

Since they are transverse to $E^u$ they must be Reebless so that
Theorem \ref{PropObienHayTorosObienTodoEsLindo} applies.

Since the foliations are $f$-invariant, the planes $P^{cs}_1$ and
$P^{cs}_2$ given by Theorem \ref{PropObienHayTorosObienTodoEsLindo} are $A$-invariant. The
fact that $P^{cs}$ contains the stable direction of $A$ is given
by Remark \ref{RemarkSubespaciosInvariantesAnosov} and Corollary
\ref{PropositionQuasiIsometriaII} since it implies that $P^{cs}$
cannot be contained in $E^u_A$.

Assume first that the planes $P^{cs}_1$ and $P^{cs}_2$ coincide.
The foliations remain at distance $R$ from translates of the
planes. By Corollary \ref{PropositionQuasiIsometriaII} we know
that two points in the same unstable leaf must separate in a
direction transverse to $P^{cs}_1=P^{cs}_2$. If $\cF^{cs}_1$ is
different from $\cF^{cs}_2$ we have a point $x$ such that
$\cF^{cs}_1(x) \neq \cF^{cs}_2(x)$. By the global product
structure we get a point $y \in \cF^{cs}_1(x)$ such that $\tilde
\cF^u(y) \cap \cF^{cs}_2(x) \neq \{y\}$. Iterating forward and
using Corollary \ref{PropositionQuasiIsometriaII} we contradict
the fact that leaves of $\cF^{cs}_1$ and $\cF_2^{cs}$ remain at
distance $R$ from translates of $P^{cs}_1=P^{cs}_2$.

Now, if $P^{cs}_1\neq P^{cs}_2$ we know that $A$ has stable
dimension $1$ since we know that $E^s_A$ is contained in both.
Using Corollary \ref{PropositionQuasiIsometriaII} and the fact
that the unstable foliation is $\tilde f$ invariant we see that
this cannot happen.

\lqqd

The following question is answered in Appendix \ref{Apendice1} for
the strong partially hyperbolic case, but in general we do not
know whether it is true:

\begin{quest}\label{QuestionPlane} In the case where $A$ has only one stable eigenvalue
(and thus, two real and different unstable eigenvalues thanks to
Proposition \ref{PropositionAnotieneinestablesComplejos}) is it
true that the plane $P^{cs}$ corresponds to the sum of the stable
eigenspace and the weak unstable one?
\end{quest}

To answer this question in Appendix \ref{Apendice1} we shall use
the following:
%
%

\begin{cor}\label{Corolario-RobustezDeDireccionTransversal}
Given a dynamically coherent partially hyperbolic diffeomorphism
$f: \TT^3 \to \TT^3$ isotopic to Anosov we know that it is
$C^1$-robustly dynamically coherent and that the
$f_\ast$-invariant plane $P$ given by Theorem
\ref{PropObienHayTorosObienTodoEsLindo} for the unique
$f$-invariant foliation $\cF^{cs}$ tangent to $E^{cs}$ does not
change for diffeomorphisms $C^1$-close to $f$.
\end{cor}

\dem The robustness of dynamical coherence follows from the fact
that being dynamically coherent it is robustly almost dynamically
coherent so that Theorem A applies.

From uniqueness and the fact that for a perturbation $g$ of $f$ we
can use the foliation $\cF^{cs}_f$ tangent to $E^{cs}_f$ as a
foliation transverse to $E^{u}_g$.  We get that the plane $P^{cs}$
almost parallel to $\cF^{cs}_f$ which is invariant under $f_\ast$
is also invariant under $g_\ast=f_\ast$. This implies that the
plane which is almost parallel to the unique $g$-invariant
foliation is again $P^{cs}$ and proves the corollary.

\lqqd


\section{Strong partially hyperbolic diffeomorphisms of $\TT^3$}\label{Section-TeoremaB}

\subsection{Strategy of the proof}
\subsubsection{} The idea of the proof of Theorem B is to obtain a global product structure
between the foliations involved in order to then get dynamical
coherence. In a certain sense, this is a similar idea to the one
used for the proof of Theorem A. However, the fact that global
product structure implies dynamical coherence is much easier in
this case due to the existence of $f$-invariant branching
foliations tangent to the center-stable direction (see subsection
\ref{SubSection-Branching}).

\subsubsection{} This approach goes in the inverse direction to the
one made in \cite{BBI2} (and continued in \cite{Hammerlindl}). In
\cite{BBI2} the proof proceeds by showing that the planes close to
the foliations are different (by using absolute domination) for
then showing (again by using absolute domination) that leaves of
$\tilde \cF^u$ are quasi-isometric so that Brin's criterium for
absolutely dominated partially hyperbolic systems (\cite{Brin})
can be applied.

Then, in \cite{Hammerlindl} it is proved that in fact, the planes
$P^{cs}$ and $P^{cu}$ close to the $f$-invariant foliations are
the expected ones in order to obtain global product structure and
then leaf conjugacy to linear models.

Another difference with their proof there is that in our case it
will be important to discuss depending on the isotopy class of
$f$. In a certain sense, the reason why in each case there is a
global product structure can be regarded as different: In the
isotopic to Anosov case (see Appendix \ref{Apendice1}) the reason
is that we deduce that the foliations are without holonomy and use
Theorem \ref{Teorema-HectorHirsch}. In the case which is isotopic
to a non-hyperbolic matrix we must find out which are the planes
close to each foliation first in order to then get the global
product structure using this fact.


\subsection{Preliminary discussions}

\subsubsection{} Let $f: \TT^3 \to \TT^3$ be a strong partially hyperbolic
diffeomorphism with splitting $T\TT^3= E^s \oplus E^c \oplus E^u$.

%

%

Corollary \ref{TeoremaMio} follows from Theorem A and the fact that strongly partially
hyperbolic diffeomorphisms are almost dynamical coherent
(\cite[Key Lemma]{BI}). We will give an
independent proof in Appendix \ref{Apendice1} since in the context
of strong partial hyperbolicity the proof becomes simpler. In this section we may safely assume then that $f_\ast$ is not Anosov.

\subsubsection{} The starting point of our proof of Theorem B is the existence of
$f$-invariant branching foliations $\cF^{cs}_{bran}$ and
$\cF^{cu}_{bran}$ tangent to $E^s\oplus E^c$ and $E^c \oplus E^u$
respectively. By using Theorem \ref{TeoBuragoIvanov} and Theorem
\ref{PropObienHayTorosObienTodoEsLindo} we can deduce the
following:

\begin{prop}\label{ProposicionDicotomiaEnLaFoliacion}
There exist an $f_\ast$-invariant plane $P^{cs}$ and $R>0$ such
that every leaf of $\tilde \cF^{cs}_{bran}$ (the lift of
$\cF^{cs}_{bran}$) lies in the $R$-neighborhood of a plane
parallel to $P^{cs}$.
Moreover, one can choose $R$ such that one of the following
conditions holds:
\begin{itemize}
\item[(i)] The projection of the plane $P^{cs}$ is dense in
$\TT^3$ and the $R$-neighborhood of every leaf of $\tilde
\cF^{cs}_{bran}$ contains a plane parallel to $P^{cs}$, or,
\item[(ii)] The projection of $P^{cs}$ is a linear two-dimensional
torus and there is a leaf of $\cF^{cs}_{bran}$ which is a
two-dimensional torus homotopic to $p(P^{cs})$.
\end{itemize}
An analogous dichotomy holds for $\cF^{cu}_{bran}$.
\end{prop}

\dem We consider sufficiently small $\eps>0$ and the foliation
$\cS_\eps$ given by Theorem \ref{TeoBuragoIvanov}.

Let $h^{cs}_\eps$ be the continuous and surjective map which is
$\eps$-close to the identity sending leaves of $\cS_\eps$ into
leaves of $\cF^{cs}_{bran}$. This implies that given a leaf $L$ of $\cF^{cs}_{bran}$ there
exists a leaf $S$ of $\cS_\eps$ such that $L$ is at distance
smaller than $L$ from $S$ and viceversa.

Since the foliation $\cS_\eps$ is transverse to $E^u$ we can apply
Theorem \ref{PropObienHayTorosObienTodoEsLindo} and we obtain that
there exists a plane $P^{cs}$ and $R>0$ such that every leaf of
the lift $\tilde \cS_\eps$ of $\cS_\eps$ to $\RR^3$ lies in an
$R$-neighborhood of a translate of $P^{cs}$.

From the previous remark, we get that every leaf of $\tilde
\cF^{cs}_{bran}$, the lift of $\cF^{cs}_{bran}$ to $\RR^3$ lies in
an $(R+ \eps)$-neighborhood of a translate of $P^{cs}$. Since
$\cF^{cs}_{bran}$ is $f$-invariant, we deduce that the plane
$P^{cs}$ is $f_\ast$-invariant.

By Proposition \ref{Proposicion-SiPesToroHayHojaToro} we know that
if $P^{cs}$ projects into a two-dimensional torus, we obtain that
the foliation $\cS_\eps$ must have a torus leaf. The image of this
leaf by $h^{cs}_\eps$ is a torus leaf of $\cF^{cs}_{bran}$.

Since a plane whose projection is not a two-dimensional torus must
be dense we get that if option (ii) does not hold, we have that
the image of $P^{cs}$ must be dense. Moreover, option (i) of
Theorem \ref{PropObienHayTorosObienTodoEsLindo} must hold for
$\cS_\eps$ and this concludes the proof of this proposition.
\lqqd

\subsection{Global product structure implies dynamical coherence}

Assume that $f: \TT^3 \to \TT^3$ is a strong partially hyperbolic
diffeomorphism. Let $\cF^{cs}_{bran}$ be the $f$-invariant
branching foliation tangent to $E^s\oplus E^c$ given by Theorem
\ref{TeoBuragoIvanov} and let $\cS_\eps$ be a foliation tangent to
an $\eps$-cone around $E^s \oplus E^c$ which remains $\eps$-close
to the lift of $\cF^{cs}_{bran}$ to the universal cover for small
$\eps$.

When the lifts of $\cS_\eps$ and $\cF^u$ to the universal cover
have a global product structure, we deduce from Proposition
\ref{PropositionQuasiIsometria} the following:

\begin{cor}\label{PropQuasiIsometria}
The foliation $\tilde \cF^u$ is quasi-isometric. Indeed, if $v \in
(P^{cs})^\perp$ is a unit vector, there exists $L>0$ such that
every unstable curve starting at a point $x$ of length larger than
$n L$  intersects $P^{cs} + n v + x$ or $P^{cs} - n v + x$.
\end{cor}

Before we show that global product structure implies dynamical
coherence, we must show that global product structure is
equivalent to a similar property with the branching foliation:

\begin{lema}\label{Lemma-GPSentreSepsybranched}
There exists $\eps>0$ such that $\tilde \cF^u$ and $\tilde
\cS_\eps$ have global product structure if and only if:
\begin{itemize}
\item[-] For every $x,y \in \RR^3$ and for every $L \in \tilde
\cF^{cs}_{bran}(y)$ we have that $\tilde \cF^u(x) \cap L \neq
\emptyset$.
\end{itemize}
\end{lema}

\dem First notice that the hypothesis imply that $\tilde \cS_\eps$
cannot have dead-end components. In particular, there exists $R>0$
and a plane $P^{cs}$ such every leaf of $\tilde \cS_\eps$ and
every leaf of $\tilde \cF^{cs}_{bran}$ verifies that it is
contained in an $R$-neighborhood of a translate of $P^{cs}$ and
the $R$-neighborhood of the leafs contains a translate of $P^{cs}$
too (see Proposition \ref{ProposicionDicotomiaEnLaFoliacion}).

We prove the direct implication first. Consider $x,y \in \RR^3$
and $L$ a leaf of $\tilde \cF^{cs}_{bran}(y)$. Now, we know that
$L$ separates in $\RR^3$ the planes $P^{cs}+ y + 2R$ and $P^{cs} +
y -2R$. One of them must be in the connected component of $\RR^3
\setminus L$ which not contains $x$, without loss of generality we
assume that it is $P^{cs} + y+ 2R$. Now, we know that there is a
leaf $S$ of $\tilde \cS_\eps$ which is contained in the half space
bounded by $P^{cs} + y +R$ not containing $L$ (notice that $L$
does not intersect $P^{cs}+ y+R$). Global product product
structure implies that $\tilde \cF^u(x)$ intersects $S$ and thus,
it also intersects $L$.

The converse direction has an analogous proof.
\lqqd

We can prove the following result which does not make use of the
isotopy class of $f$.

\begin{prop}\label{Proposition-GPSIMPLICACOHERENCIA}
Assume that there is a global product structure between the lift
of $\cS_\eps$ and the lift of $\cF^u$ to the universal cover. Then
there exists an $f$-invariant foliation $\cF^{cs}$ everywhere
tangent to $E^s \oplus E^u$.
\end{prop}

\dem We will show that the branched foliation $\tilde
\cF^{cs}_{bran}$ must be a true foliation. Using Proposition
\ref{PropositionBranchingSinBranchingEsFoliacion} this is reduced
to showing that each point in $\RR^3$ belongs to a unique leaf of
$\tilde \cF^{cs}_{bran}$.

Assume otherwise, i.e. there exists $x\in \RR^3$ such that $\tilde
\cF^{cs}_{bran}(x)$ has more than one complete surface. We call
$L_1$ and $L_2$ different leaves in $\tilde \cF^{cs}_{bran}(x)$.
There exists $y$ such that $y \in L_1 \setminus L_2$. Using global
product structure and Lemma \ref{Lemma-GPSentreSepsybranched} we
get $z \in L_2$ such that:

\begin{itemize}
\item[-] $y \in \tilde \cF^u(z)$.
\end{itemize}

Consider $\gamma$ the arc in $\tilde \cF^u(z)$ whose endpoints are
$y$ and $z$. Let $R$ be the value given by Proposition
\ref{ProposicionDicotomiaEnLaFoliacion} and $\ell>0$ given by
Corollary \ref{PropQuasiIsometria}. We consider $N$ large enough
so that $\tilde f^N(\gamma)$ has length larger than $n \ell$ with
$n \gg R$.

By Corollary \ref{PropQuasiIsometria} we get that the distance
between $P^{cs}+ \tilde f^N(z)$ and $\tilde f^N(y)$ is much larger
than $R$. However, we have that, by $\tilde f$-invariance of
$\tilde \cF^{cs}_{bran}$ there is a leaf of $\tilde
\cF^{cs}_{bran}$ containing both $\tilde f^N(z)$ and $\tilde
f^N(x)$ and another one containing both $\tilde f^N(y)$ and
$\tilde f^N(x)$. This contradicts Proposition
\ref{ProposicionDicotomiaEnLaFoliacion} showing that $\tilde
\cF^{cs}_{bran}$ must be a true foliation.

\lqqd

\subsection{Torus leafs.}

\subsubsection{} This subsection is devoted to prove the following:

\begin{lema}\label{LemaNohayToros}
If $\cF^{cs}_{bran}$ contains a leaf which is a two-dimensional
torus, then there is a leaf of $\cF^{cs}_{bran}$ which is a torus
and it is fixed by $f^k$ for some $k$. Moreover, this leaf is
repelling.
\end{lema}

\dem Let $T \en \TT^3$ be a leaf of $\cF^{cs}$ homeomorphic to a
two-torus. Since $\cF^{cs}$ is $f$-invariant and $P^{cs}$ is
invariant under $f_\ast$ we get that the image of $T$ by $f$ is
homotopic to $T$ and a leaf of $\cF^{cs}_{bran}$.

Notice that having an $f_\ast$-invariant plane which projects into
a torus already implies that $f_\ast$-cannot be hyperbolic (see
Proposition \ref{PropAnosovSonIrreducibles}).

By Lemma \ref{LemaPlanosInvariantes}, the plane $P^{cs}$ coincides
with $E^s_\ast \oplus E^u_\ast$ (the eigenspaces corresponding to
the eigenvalues of modulus different from $1$ of $f_\ast$)

Since the eigenvalue of $f_\ast$ in $E^c_\ast$ is of modulus $1$,
this implies that if we consider two different lifts of $T$, then
they remain at bounded distance when iterated by $\tilde f$.
Indeed, if we consider two different lifts $\tilde T_1$ and
$\tilde T_2$ of $T$ we have that $\tilde T_2 = \tilde T_1 +\gamma$
with $\gamma \in E^c_\ast \cap \ZZ^3$. Now, we have that $\tilde f
(\tilde T_2) = \tilde f(\tilde T_1) + f_\ast(\gamma) = \tilde
f(\tilde T_1) \pm \gamma$.

We shall separate the proof depending on how the orbit of $T$ is.
\medskip

{\bf Case 1:} Assume the torus $T$ is fixed by some iterate $f^n$
of $f$. Then, since it is tangent to the center stable
distribution, we obtain that it must be repelling as desired.
\medskip

{\bf Case 2:} If the orbit of $T$ is dense, we get that
$\cF^{cs}_{bran}$ is a true foliation by two-dimensional tori
which we call $\cF^{cs}$ from now on. This is obtained by the fact
that one can extend the foliation to the closure using the fact
that there are no topological crossings between the torus leaves
(see Proposition \ref{Proposition-BWProp16}).

Since all leaves must be two-dimensional torus  which are
homotopic we get that the foliation $\cF^{cs}$ has no holonomy see
Proposition \ref{Proposition-ReebStability} and Proposition
\ref{Proposition-SucesionDeTorosEnBranchingFol}).

Using Theorem \ref{Teorema-HectorHirsch}, we get that the unstable
direction $\tilde \cF^u$ in the universal cover must have a global
product structure with $\tilde \cF^{cs}$.

Let $S$ be a leaf of $\cF^{cs}$ and consider $\tilde S_1$ and
$\tilde S_2$ two different lifts of $S$ to $\RR^3$.

Consider an arc $J$ of $\tilde \cF^u$ joining $\tilde S_1$ to
$\tilde S_2$. Iterating the arc $J$ by $\tilde f^n$ we get that
its length grows exponentially, while the extremes remain the the
forward iterates of $\tilde S_1$ and $\tilde S_2$ which remain at
bounded distance by the argument above.

By considering translations of one end of $\tilde f^n(J)$ to a
fundamental domain and taking a convergent subsequence we obtain a
leaf of $\tilde \cF^u$ which does not intersect every leaf of
$\tilde \cF^{cs}$. This contradicts global product structure.

\medskip

{\bf Case 3:} Let $T_1, T_2 \in \cF^{cs}_{bran}$ two different
torus leaves. Since there are no topological crossings, we can
regard $T_2$ as embedded in $\TT^2 \times [-1,1]$ where both
boundary components are identified with $T_1$ and such that the
embedding is homotopic to the boundary components (recall that any
pair of torus leaves must be homotopic). In particular, we get
that $\TT^3 \setminus (T_1 \cup T_2)$ has at least two different
connected components and each of the components has its boundary
contained in $T_1 \cup T_2$.

If the orbit of $T$ is not dense, we consider
$\cO=\overline{\bigcup_n f^n(T)}$ the closure of the orbit of $T$
which is an invariant set.

Recall that we can assume completeness of $\cF^{cs}_{bran}$ (i.e.
for every $x_n \to x$ and $L_n \in \cF^{cs}_{bran}(x_n)$ we have
that $L_n$ converges in the $C^1$-topology to $L_\infty \in
\cF^{cs}_{bran}(x)$). We get that $\cO$ is saturated by leaves of
$\cF^{cs}_{bran}$ all of which are homotopic torus leaves (see
Proposition \ref{Proposition-SucesionDeTorosEnBranchingFol}).

Let $U$ be a connected component of the complement of $\cO$. By
the previous remarks we know that its boundary $\partial U$ is
contained in the union of two torus leaves of $\cF^{cs}_{bran}$.

If some component $U$ of $\cO^c$ verifies that there exists $n\geq
1$ such that $f^n(U) \cap U \neq \emptyset$, by invariance of
$\cO^c$ we get that $f^{2n}$ fixes both torus leaves whose union
contains $\partial U$. This implies the existence of a periodic
normally repelling torus as in Case 1.

We claim that if every connected component of $\cO^c$ is
wandering, then we can show that every leaf of $\tilde \cF^u$
intersects every leaf of $\tilde \cF^{cs}_{bran}$ which allows to
conclude exactly as in Case 2.

To prove the claim, consider $\delta$ given by the local product
structure between these two transverse foliations (one of them
branched). This means that given $x,y$ such that $d(x,y)<\delta$
we have that $\tilde \cF^u(x)$ intersects every leaf of $\tilde
\cF^{cs}_{bran}$ passing through $y$.

Assume there is a point $x\in \RR^3$ such that $\tilde \cF^u(x)$
does not intersect every leaf of $\tilde \cF^{cs}_{bran}$. We know
that each leaf of $\tilde \cF^{cs}_{bran}$ separates $\RR^3$ into
two connected components so we can choose among the lifts of torus
leaves, the leaf $\tilde T_0$ which is the lowest (or highest
depending on the orientation of the semi-unstable leaf of $x$ not
intersecting every leaf of $\tilde \cF^{cs}_{bran}$) not
intersecting $\tilde \cF^{u}(x)$. We claim that $\tilde T_0$ must
project by the covering projection into a torus leaf which
intersects\footnote{If it were a true foliation, the leaf would
coincide with the boundary. In this case, it could be that the
boundary of the connected component of $\cO^c$ is only a part of
the torus. In any case, this is enough for our argument.} the
boundary of a connected component of $\cO^c$. Indeed, there are
only finitely many connected components $U_1, \ldots, U_N$ of
$\cO^c$ having volume smaller than the volume of a $\delta$-ball,
so if a point is not in $U_i$ for some $i$, we know that it must
be covered by local product structure boxes forcing its unstable
leaf to advance until one of those components.

On the other hand, using $f$-invariance of $\cF^{u}$ and the fact
that every connected component of $\cO^c$ is wandering, we get
that every point in $U_i$ must eventually fall out of $\bigcup_i
U_i$ and then its unstable manifold must advance to other
component. This concludes the claim, and as we explained, allows
to use the same argument as in Case 2 to finish the proof in Case
3. \lqqd

M.A. Rodriguez Hertz and R. Ures were kind to communicate an
alternative proof of this lemma by using an adaptation of an
argument due to Haefliger for branching foliations (it should
appear in \cite{HHU}).

\subsection{Existence of a Global Product Structure}

\subsubsection{}
In this section we will prove the following result which will
allow us to conclude in the case where $f_\ast$ is not isotopic to
Anosov.

\begin{prop}\label{PropGPSyPlanos}
Let $f: \TT^3 \to \TT^3$ be a strongly partially hyperbolic
diffeomorphism which is not isotopic to Anosov. Assume moreover
that $\cF^{cs}_{bran}$ does not contain a torus leaf. Then, the
plane $P^{cs}$ given by Proposition
\ref{ProposicionDicotomiaEnLaFoliacion} corresponds to the
eigenplane corresponding to the eigenvalues of modulus smaller or
equal to $1$. Moreover, there is a global product structure
between $\tilde \cF^{cs}_{bran}$ and $\tilde \cF^u$. A symmetric
statement holds for $\tilde \cF^{cu}_{bran}$ and $\tilde \cF^s$.
\end{prop}

\subsubsection{} Since $f_\ast$ is not Anosov, we have from subsection \ref{SubSection-ClaseIsotopPH} that there
are exactly three $f_\ast$-invariant lines $E^s_\ast$, $E^c_\ast$
and $E^u_\ast$ corresponding to the eigenvalues of $f_\ast$ of
modulus smaller, equal and larger than one respectively.


\begin{lema}\label{LemaEstableNoSeQuedaCercaDeCentroInestable}
For every $R>0$ and $x\in \RR^3$ we have that $\tilde \cF^u(x)$ is
not contained in an $R$-neighborhood of $(E^s_\ast \oplus
E^c_\ast) + x$. Symmetrically, for every $R>0$ and $x \in \RR^3$
the leaf $\tilde \cF^s(x)$ is not contained in an $R$-neighborhood
of $(E^c_\ast \oplus E^u_\ast) +x$.
\end{lema}

\dem Let $C$ be a connected set contained in an $R$-neighborhood
of a translate of $E^s_\ast \oplus E^c_\ast$, we will estimate the
diameter of $\tilde f(C)$ in terms of the diameter of $C$.

\begin{af}
There exists $K_R$ which depends only on $\tilde f$, $f_\ast$ and
$R$ such that:
$$ \diam (\tilde f(C)) \leq \diam(C) + K_R $$
\end{af}

\dem  Let $K_0$ be the $C^0$-distance between $\tilde f$ and
$f_\ast$ and consider $x,y \in C$ we get that:

$$ d(\tilde f(x),\tilde f(y)) \leq d(f_\ast(x),f_\ast(y)) + d(f_\ast(x),\tilde f(x))+ d(f_\ast(y),\tilde f(y)) \leq $$
$$ \leq d(f_\ast(x), f_\ast(y)) + 2 K_0$$

We have that the difference between $x$ and $y$ in the unstable
direction of $f_\ast$ is bounded by $2R$ given by the distance to
the plane $E^s_\ast \oplus E^u_\ast$ which is transverse to
$E^u_\ast$.

We then have that $d(f_\ast(x),f_\ast(y)) \leq d(x,y) + 2
|\lambda^u| R$ where $\lambda^u$ is the eigenvalue of modulus
larger than $1$ obtaining:

$$ d(\tilde f(x),\tilde f(y))  \leq d(x,y) + 2K_0 + 2|\lambda^u| R = d(x,y) + K_R$$

\noindent which concludes the proof of the claim. \finobs

Now, this implies that if we consider an arc $\gamma$ of $\tilde
\cF^u$ of length $1$ and assume that its future iterates remain in
a slice parallel to $E^s_\ast \oplus E^c_\ast$ of width $2R$ we
have that

$$\diam( \tilde f^n(\gamma)) < \diam(\gamma) + n K_R \leq 1 + n K_R$$

So that the diameter grows linearly with $n$.

The volume of balls in the universal cover of $\TT^3$ grows
polynomially with the radius (see Step 2 of \cite{BBI} or page 545
of \cite{BI}, notice that the universal) so that we have that
$B_\delta(\tilde f^{n}(\gamma))$ has volume which is polynomial
$P(n)$ in $n$.

On the other hand, we know from the partial hyperbolicity that
there exists $C>0$ and $\lambda>1$ such that the length of $\tilde
f^n(\gamma)$ is larger than $C \lambda^n$.

Using Theorem \ref{CorolarioConsecuenciasReeb} (iv), we obtain
that there exists $n_0$ uniform such that every arc of length $1$
verifies that $\tilde f^{n_0}(\gamma)$ is not contained in the
$R$-neighborhood of a translate of $E^s_\ast \oplus E^c_\ast$.
This implies that no unstable leaf can be contained in the
$R$-neighborhood of a translate of $E^s_\ast \oplus E^c_\ast$
concluding the proof of the lemma.

\lqqd

\subsubsection{} We are now ready to prove Proposition \ref{PropGPSyPlanos}

\demo{of Proposition \ref{PropGPSyPlanos}} Consider the plane
$P^{cs}$ given by Proposition
\ref{ProposicionDicotomiaEnLaFoliacion} for the branching
foliation $\cF^{cs}_{bran}$.

If option (ii) of Proposition
\ref{ProposicionDicotomiaEnLaFoliacion} holds, we get that there
must be a torus leaf in $\cF^{cs}_{bran}$ which we assume there is
not.

By Lemma \ref{LemaPlanosInvariantes} and Proposition
\ref{ProposicionDicotomiaEnLaFoliacion} the plane $P^{cs}$ must be
either $E^s_\ast \oplus E^c_\ast$ or $E^c_\ast \oplus E^u_\ast$.

Lemma \ref{LemaEstableNoSeQuedaCercaDeCentroInestable} implies
that $P^{cs}$ cannot be $E^c_\ast \oplus E^u_\ast$ since $\cF^s$
is contained in $\cF^{cs}_{bran}$. This implies that $P^{cs}=
E^s_\ast \oplus E^c_\ast$ as desired.

Now, using Lemma \ref{LemaEstableNoSeQuedaCercaDeCentroInestable}
for $\tilde \cF^u$ we see that the unstable foliation cannot
remain close to a translate of $P^{cs}$ and must intersect every
leaf of $\tilde \cF^{cs}_{bran}$ obtaining the desired global
product structure.

\lqqd

\subsection{Proof of Theorem B}\label{SectionCoherenciaPrueba}

To prove Theorem B, we first assume that $f_\ast$ is not isotopic
to Anosov.

Consider the branching foliation $\cF^{cs}_{bran}$ given by
Theorem \ref{TeoBuragoIvanov}. We can apply Proposition
\ref{ProposicionDicotomiaEnLaFoliacion} to $\cF^{cs}_{bran}$ and
obtain a plane $P^{cs}$ which is close to the lift of leaves of
$\cF^{cs}_{bran}$.

If the plane $P^{cs}$ projects into a torus, there must be a
two-dimensional torus as a leaf of $\cF^{cs}_{bran}$ then, by
Lemma \ref{LemaNohayToros} we obtain a repelling periodic
two-dimensional torus.

If $P^{cs}$ is not a torus, then Proposition \ref{PropGPSyPlanos}
applies giving a global product structure between the lift of the
unstable foliation and the lift of $\cF^{cs}_{bran}$.

By Proposition \ref{Proposition-GPSIMPLICACOHERENCIA} we get that
there exists an $f$-invariant foliation $\cF^{cs}$ tangent to $E^s
\oplus E^c$.

The proof shows that there must be a unique $f$-invariant
foliation tangent to $E^{cs}$ (and to $E^{cu}$).

Indeed, we get that every foliation tangent to $E^{cs}$ must
verify option (i) of Proposition
\ref{ProposicionDicotomiaEnLaFoliacion} when lifted to the
universal cover and that the plane which is close to the foliation
must correspond to the eigenspace of $f_\ast$ corresponding to the
smallest eigenvalues (Proposition \ref{PropGPSyPlanos}).

Using quasi-isometry of the strong foliations, this implies that
if there is another surface tangent to $E^{cs}$ through a point
$x$, then this surface will not extend to an $f$-invariant
foliation since we get that forward iterates will get arbitrarily
far from this plane.

This concludes the proof of Theorem B in case $f$ is not isotopic
to Anosov. The other case is covered by Theorem A (see also the
appendix for a simpler proof in the strong partially hyperbolic
setting).

\lqqd


\appendix

\section{Strong partially hyperbolic diffeomorphisms in the isotopy class of Anosov.}\label{Apendice1}

In this appendix we prove the following:

\begin{teo}\label{TeoremaMio} Let $f:\TT^3 \to \TT^3$ a strong partially hyperbolic diffeomorphism isotopic to a
linear Anosov automorphism. Then, $f$ is dynamically coherent.
\end{teo}

This result can be seen as a direct consequence of Theorem A and
the fact that every strong partially hyperbolic diffeomorphism is
almost dynamically coherent (\cite{BI}).

We present here a simpler proof of this result which is
independent of Section \ref{Section-PHAnosov}. Then, we prove a
result in the vein of Proposition \ref{PropGPSyPlanos} for the
isotopy class of Anosov (the difference is in this case it is an a
fortiori result while in the other case it is needed to get
dynamical coherence), the proof of this last result is based on
what is proved in Section \ref{Section-PHAnosov} and the main
result from \cite{DW}.

\demo{of Theorem \ref{TeoremaMio}} Notice that if $f_\ast$ is
hyperbolic, then, every invariant plane must be totally irrational
(see Remark \ref{RemarkSubespaciosInvariantesAnosov}), so that it
projects into a plane in $\TT^3$.

Let $\cF^{cs}_{bran}$ be the branched foliation tangent to
$E^{cs}$ given by Theorem \ref{TeoBuragoIvanov}. By Proposition
\ref{ProposicionDicotomiaEnLaFoliacion} we get a
$f_\ast$-invariant plane $P^{cs}$ in $\RR^3$ which we know cannot
project into a two-dimensional torus since $f_\ast$ has no
invariant planes projecting into a torus, this implies that option
(i) of Proposition \ref{ProposicionDicotomiaEnLaFoliacion} is
verified.

Since for every $\eps>0$, Theorem \ref{TeoBuragoIvanov} gives us a
foliation $\cS_\eps$ whose lift is close to $\tilde \cF^{cs}$, we
get that the foliation $\tilde \cS_\eps$ remains close to $P^{cs}$
which must be totally irrational. By Lemma
\ref{RemarkHojasCerradas} (i) we get that all leaves of $\cS_\eps$
are simply connected, thus, we get that the foliation $\cS_\eps$
is without holonomy.

We can apply Theorem \ref{Teorema-HectorHirsch} and we obtain that
for every $\eps>0$ there is a global product structure between
$\tilde \cS_\eps$ and $\tilde \cF^u$ which is transverse to
$\cS_\eps$ if $\eps$ is small enough.

The rest of the proof follows from Proposition
\ref{Proposition-GPSIMPLICACOHERENCIA}. \lqqd

In fact, using the same argument as in Proposition
\ref{Prop-UnicidadDeFoliacionTgAEcs} we get uniqueness of the
foliation tangent to $E^s\oplus E^c$.

We are also able to prove the following proposition which is
similar to Proposition \ref{PropGPSyPlanos} in the context of
partially hyperbolic diffeomorphisms isotopic to Anosov. This will
be used in \cite{HP} to obtain leaf conjugacy to the linear model.

Notice first that the eigenvalues of $f_\ast$ verify that they are
all different (see Propositions \ref{PropAnosovSonIrreducibles}
and \ref{PropositionAnotieneinestablesComplejos}).

We shall name them $\lambda_1, \lambda_2, \lambda_3$ and assume
they verify:

$$ |\lambda_1|< |\lambda_2|< |\lambda_3| \quad ; \quad |\lambda_1|<1 \ , \ |\lambda_2|\neq 1 \ , \ |\lambda_3|>1 $$

\noindent we shall denote as $E^i_\ast$ to the eigenline of
$f_\ast$ corresponding to $\lambda_i$.

\begin{prop}\label{PropPlanoCasoIsotopico}
The plane close to the branched foliation $\tilde \cF^{cs}$
corresponds to the eigenplane corresponding to the eigenvalues of
smaller modulus (i.e. the eigenspace $E^{1}_\ast \oplus E^2_\ast$
corresponding to $\lambda_1$ and $\lambda_2$). Moreover, there is
a global product structure between $\tilde \cF^{cs}$ and $\tilde
\cF^u$. A symmetric statement holds for $\tilde \cF^{cu}$ and
$\tilde \cF^s$.
\end{prop}

\dem  This proposition follows from the existence of a
semiconjugacy $H$ between $\tilde f$ and its linear part $f_\ast$
which is at bounded distance from the identity.

The existence of a global product structure was proven above.
Assume first that $|\lambda_2|<1$, in this case, we know that
$\tilde \cF^u$ is sent by the semiconjugacy into lines parallel to
the eigenspace of $\lambda_3$ for $f_\ast$. This readily implies
that $P^{cs}$ must coincide with the eigenspace of $f_\ast$
corresponding to $\lambda_1$ and $\lambda_2$ otherwise we would
contradict the global product structure.

The case were $|\lambda_2|>1$ is more difficult. First, it is not
hard to show that the eigenspace corresponding to $\lambda_1$ must
be contained in $P^{cs}$ (otherwise we can repeat the argument in
Lemma \ref{LemaEstableNoSeQuedaCercaDeCentroInestable} to reach a
contradiction).

Assume by contradiction that $P^{cs}$ is the eigenspace
corresponding to $\lambda_1$ and $\lambda_3$.

First, notice that by the basic properties of the semiconjugacy
$H$, for every $x\in \RR^3$ we have that $\tilde \cF^u(x)$ is sent
by $H$ into $E^u_\ast + H(x)$ (where $E^u_\ast = E^2_\ast \oplus
E^3_\ast$ is the eigenspace corresponding to $\lambda_2$ and
$\lambda_3$ of $f_\ast$).

We claim that this implies that in fact $H(\tilde \cF^u(x)) =
E^2_\ast + H(x)$ for every $x\in \RR^3$. In fact, we know from
Corollary \ref{PropositionQuasiIsometriaII} that points of
$H(\tilde \cF(x))$ which are sufficiently far apart are contained
in a cone of $(E^2_\ast \oplus E^3_\ast) + H(x)$ bounded by two
lines $L_1$ and $L_2$ which are transverse to $P^{cs}$. If
$P^{cs}$ contains $E^3_\ast$ this implies that if one considers
points in the same unstable leaf which are sufficiently far apart,
then their image by $H$ makes an angle with $E^3_\ast$ which is
uniformly bounded from below. If there is a point $y \in \tilde
\cF^u(x)$ such that $H(y)$ not contained in $E^2_\ast$ then we
have that $d(\tilde f^n(y), \tilde f^n(x))$ goes to $\infty$ with
$n$ while the angle of $H(y) - H(x)$ with $E^3_\ast$ converges to
$0$ exponentially contradicting Corollary
\ref{PropositionQuasiIsometriaII}.

Consider now a point $x\in \RR^3$ and let $y$ be a point which can
be joined to $x$ by a finite set of segments $\gamma_1, \ldots,
\gamma_k$ tangent either to $E^s$ or to $E^u$ (an $su$-path, see
\cite{DW}). We know that each $\gamma_i$ verifies that
$H(\gamma_i)$ is contained either in a translate of $E^1_\ast$
(when $\gamma_i$ is tangent to $E^s$, i.e. it is an arc of the
strong stable foliation $\tilde \cF^s$) or in a translate of
$E^2_\ast$ (when $\gamma_i$ is tangent to $E^u$ from what we have
shown in the previous paragraph). This implies that the
\emph{accessibility} class of $x$ (see \cite{DW,HHUAcc} for a
definition and properties) verifies that its image by $H$ is
contained in $(E^1_\ast \oplus E^2_\ast) + H(x)$. The projection
of $E^1_\ast \oplus E^2_\ast$ to the torus is not the whole
$\TT^3$ so in particular, we get that $f$ cannot be accessible.
From Corollary \ref{Corolario-RobustezDeDireccionTransversal} this
situation should be robust under $C^1$-perturbations since those
perturbations cannot change the direction of $P^{cs}$.

On the other hand, in \cite{DW,HHUAcc} it is proved that by an
arbitrarily small ($C^1$ or $C^r$) perturbation of $f$ one can
make it accessible.This gives a contradiction and concludes the
proof.

\lqqd




\end{document}

%% file: fmasopciones.pstex_t
\begin{picture}(0,0)%
\includegraphics{fmasopciones.pstex}%
\end{picture}%
\setlength{\unitlength}{3947sp}%
\begingroup\makeatletter\ifx\SetFigFont\undefined%
\gdef\SetFigFont#1#2#3#4#5{%
  \reset@font\fontsize{#1}{#2pt}%
  \fontfamily{#3}\fontseries{#4}\fontshape{#5}%
  \selectfont}%
\fi\endgroup%
\begin{picture}(6529,2247)(264,-1560)
\put(2339,-1499){\makebox(0,0)[lb]{\smash{{\SetFigFont{12}{14.4}{\rmdefault}{\mddefault}{\updefault}$x$}}}}
\put(4801,-1161){\makebox(0,0)[lb]{\smash{{\SetFigFont{12}{14.4}{\rmdefault}{\mddefault}{\updefault}$x+\gamma$}}}}
\put(1451,-361){\makebox(0,0)[lb]{\smash{{\SetFigFont{12}{14.4}{\rmdefault}{\mddefault}{\updefault}$F_+(x)$}}}}
\put(5401,-49){\makebox(0,0)[lb]{\smash{{\SetFigFont{12}{14.4}{\rmdefault}{\mddefault}{\updefault}$F_+(x)+\gamma$}}}}
\end{picture}%

%% file: foliations.pstex_t
\begin{picture}(0,0)%
\includegraphics{foliations.pstex}%
\end{picture}%
\setlength{\unitlength}{3947sp}%
\begingroup\makeatletter\ifx\SetFigFont\undefined%
\gdef\SetFigFont#1#2#3#4#5{%
  \reset@font\fontsize{#1}{#2pt}%
  \fontfamily{#3}\fontseries{#4}\fontshape{#5}%
  \selectfont}%
\fi\endgroup%
\begin{picture}(6753,2611)(339,-2273)
\end{picture}%

%% file: curvaeta.pstex_t
\begin{picture}(0,0)%
\includegraphics{curvaeta.pstex}%
\end{picture}%
\setlength{\unitlength}{3947sp}%
\begingroup\makeatletter\ifx\SetFigFont\undefined%
\gdef\SetFigFont#1#2#3#4#5{%
  \reset@font\fontsize{#1}{#2pt}%
  \fontfamily{#3}\fontseries{#4}\fontshape{#5}%
  \selectfont}%
\fi\endgroup%
\begin{picture}(4028,2586)(1425,-2396)
\put(1764,-2236){\makebox(0,0)[b]{\smash{{\SetFigFont{12}{14.4}{\rmdefault}{\mddefault}{\updefault}$z$}}}}
\put(4101,-2174){\makebox(0,0)[b]{\smash{{\SetFigFont{12}{14.4}{\rmdefault}{\mddefault}{\updefault}$\eta_z$}}}}
\put(3126,-649){\makebox(0,0)[b]{\smash{{\SetFigFont{12}{14.4}{\rmdefault}{\mddefault}{\updefault}$z_2$}}}}
\put(2726,-199){\makebox(0,0)[b]{\smash{{\SetFigFont{12}{14.4}{\rmdefault}{\mddefault}{\updefault}$z_1$}}}}
\end{picture}%

%% file: cono.pstex_t
\begin{picture}(0,0)%
\includegraphics{cono.pstex}%
\end{picture}%
\setlength{\unitlength}{3947sp}%
\begingroup\makeatletter\ifx\SetFigFont\undefined%
\gdef\SetFigFont#1#2#3#4#5{%
  \reset@font\fontsize{#1}{#2pt}%
  \fontfamily{#3}\fontseries{#4}\fontshape{#5}%
  \selectfont}%
\fi\endgroup%
\begin{picture}(6735,4626)(489,-3777)
\put(5026,-3235){\makebox(0,0)[lb]{\smash{{\SetFigFont{12}{14.4}{\rmdefault}{\mddefault}{\updefault}$\tilde \cF^{cs}(x)$}}}}
\put(3614,-2448){\makebox(0,0)[lb]{\smash{{\SetFigFont{12}{14.4}{\rmdefault}{\mddefault}{\updefault}$x$}}}}
\put(5326,-535){\makebox(0,0)[lb]{\smash{{\SetFigFont{12}{14.4}{\rmdefault}{\mddefault}{\updefault}$L_1$}}}}
\put(2801,-160){\makebox(0,0)[lb]{\smash{{\SetFigFont{12}{14.4}{\rmdefault}{\mddefault}{\updefault}$L_2$}}}}
\put(3739,689){\makebox(0,0)[lb]{\smash{{\SetFigFont{12}{14.4}{\rmdefault}{\mddefault}{\updefault}$\tilde \cF^u(x)$}}}}
\put(6276,-361){\makebox(0,0)[lb]{\smash{{\SetFigFont{12}{14.4}{\rmdefault}{\mddefault}{\updefault}$E^u_A + x$}}}}
\end{picture}%